\documentclass{amsart}
\usepackage{amssymb}
\usepackage{amsthm}
\usepackage{hyperref}
\usepackage{amsmath}
\usepackage{amssymb}
\usepackage{enumerate}
\usepackage{graphicx}
\usepackage{amsthm}
\usepackage{amsopn}
\usepackage{amsfonts}
\usepackage{upgreek}
\usepackage{amsfonts}
\usepackage{amssymb}
\usepackage{empheq}
\usepackage{caption}
\numberwithin{equation}{section}

\newtheorem{theorem}{Theorem}
\newtheorem{corollary}{Corollary}
\newtheorem{example}{Example}%
\newtheorem{remark}{Remark}%
\newtheorem{lemma}{Lemma}%
\newtheorem{definition}{Definition}%
\begin{document}

\title[Floquet Theory for nonautonomous linear periodic IDEPCAG]{A Floquet-Lyapunov Theory for nonautonomous linear periodic differential equations with piecewise constant deviating arguments}
\author{Ricardo Torres}
\address{Instituto de Ciencias Físicas y Matemáticas, Facultad de Ciencias, Universidad Austral de Chile\\
Campus Isla Teja s/n, Valdivia, Chile.}
\address{\noindent Instituto de Ciencias, Universidad Nacional de General Sarmiento\\
Los Polvorines, Buenos Aires, Argentina}
\curraddr{}
\email{ricardo.torres@uach.cl}
\thanks{}
\subjclass[2020]{34A36, 34K11,  34K45, 39A06, 39A21}

\keywords{Piecewise constant argument, linear functional differential equations, Floquet theorem, Impulsive Differential equations, Hybrid dynamics, Periodic systems, Floquet-Lyapunov transformation.}

\date{}

\dedicatory{Dedicated to the memory of Prof. Istv\'an Gy\H ori and Prof. Nicol\'as Yus Su\'arez.}

\begin{abstract}
We present a version of the classical Floquet-Lyapunov theorem for $\omega-$periodic nonautonomous linear (impulsive and non-impulsive) differential equations with piecewise constant arguments of generalized type (in short, IDEPCAG or DEPCAG). We have proven that the nonautonomous linear IDEPCAG is kinematically similar to an autonomous linear ordinary differential equation. We have also provided some examples to demonstrate the effectiveness of our results.
\end{abstract}
\maketitle

\section{Introduction}
Discontinuous phenomena are often in nature, and they need to be represented with piecewise constant functions and impulses to illustrate an abrupt change in the state of the phenomena in study. Differential equations with deviating arguments, such as $f(t)=[t+1],$ (the greatest integer function), were analyzed by A. Myshkis in \cite{203} (1977). An example of such an equation corresponds to $$x'(t)=f(t,x(t),x([t+1])).$$

M. Akhmet proposed a generalized form of differential equations with step functions as deviating arguments in the form of 
\begin{equation}
z^{\prime}(t)=f(t,z(t),z(\gamma(t))),\label{depcag_eq}
\end{equation}
where $\gamma(t)$ is a \emph{piecewise constant argument of generalized type}.

Consider sequences $\left(t_{n}\right)_{n\in\mathbb{Z}}$ and $\left(\zeta_{n}\right)_{n\in\mathbb{Z}}$ such that $t_{n}<t_{n+1}$ for all $n\in\mathbb{Z}$, and $\displaystyle{\lim_{n\rightarrow\pm\infty}t_{n}=\pm\infty}$, with $\zeta_{n}\in[t_{n},t_{n+1}]$. Define $\gamma(t)=\zeta_{n}$ if $t\in I_{n}=\left[t_{n},t_{n+1}\right)$. In other words, $\gamma(t)$ is a step function, for example, $\gamma(t) =[t]$, where $[\cdot]$ denotes the greatest integer function, which is constant in every interval $[n,n+1[$ with $n\in \mathbb{Z}$ (see \eqref{idepca_parte_entera}).

If a $\gamma$ function is used, the interval $I_n$ is decomposed into advanced and retarded subintervals $I_{n}=I_{n}^{+}\bigcup I_{n}^{-}$, where $I_{n}^{+}=[t_{n},\zeta_{n}]$ and $I_{n}^{-}=[\zeta_{n},t_{n+1}]$. This type of differential equation is called Differential Equations with Piecewise Constant Argument of Generalized Type (DEPCAG). They have remarkable properties, as the solutions remain continuous functions, even when $\gamma$ is discontinuous. We can define a difference equation by assuming continuity of the solutions of \eqref{depcag_eq} and integrating from $t_n$ to $t_{n+1}$. Therefore, this type of differential equation has hybrid dynamics (see \cite{AK2, P2011, Torres_Castillo_Pinto_2023}).\\

If an impulsive condition is considered at instants 
$\{t_n\}_{n\in\mathbb{Z}}$, we define the 
\emph{Impulsive differential equations with piecewise constant argument of generalized type} (\emph{IDEPCAG}) (see \cite{AK3}),
\begin{align}
&z^{\prime}(t)=f(t,z(t),z(\gamma(t))),\qquad \qquad \qquad \quad t\neq t_{n} \nonumber \\
&\Delta z(t_{n})\coloneqq z(t_{n})-z(t_{n}^{-})=J_{n}(z(t_{n}^{-})),\qquad t=t_{n},\quad n\in\mathbb{N} \label{idepcag_gral},
\end{align}
where $z(t_n^-)=\displaystyle{\lim_{t\to t_n^-}z(t),}$ and $J_n$ is the impulsive operator (see \cite{Samoilenko}).\\
When the differential equation explicitly shows the piecewise constant argument used, we will call it DEPCA (or IDEPCA if it has impulses).\\

Let the following ordinary differential system
\begin{align}
x^{\prime}(t)=A(t)x(t),\quad A(t+\omega)=A(t),\quad \forall t\in\mathbb{R},\label{periodico_ordinario}
\end{align}
where $A(t)$ is a continuous matrix. What can be said about the stability of solutions?
The following example demonstrates that the eigenvalues are insufficient to ensure solution stability:
\begin{example}{(Counterexample of \emph{Markus-Yamabe})}\cite{markus_yamabe}\\
Let the system
\begin{eqnarray}
 x'=A(t)x,\quad A(t+\pi)=A(t),\label{markus_yamabe}
\end{eqnarray}
where 
$$A(t)=
\begin{pmatrix}
-1+\frac{3}{2}\cos^2(t) & 1-\frac{3}{2}\sin(t)\cos(t)\\
&\\
-1-\frac{3}{2}\sin(t)\cos(t) & -1+\frac{3}{2}\sin^2(t)
\end{pmatrix}
.$$
The matrix $A(t)$ has eigenvalues that are constant and equal to $\frac{1}{4}\left(-1\pm\sqrt{7}i\right).$ At first glance, we might conclude that the zero solution of equation \eqref{markus_yamabe} is asymptotically stable due to the negative real part of the eigenvalues. However, a solution of the same equation is given by 
$$x(t)=\exp{(t/2)}
\begin{pmatrix}
-\cos(t) & \\
\quad\sin(t) &
\end{pmatrix}
,$$
which is unbounded. Therefore, the zero solution of \eqref{markus_yamabe} is unstable.
\end{example}
\noindent Consequently, a natural question arises:
\begin{center}
¿What can be said about the stability of a nonautonomous linear system using its eigenvalues?    
\end{center}
In an attempt to study the stability of \eqref{periodico_ordinario} with the classical autonomous spectral theory, the French mathematician \textit{G. Floquet} proved, in 1883, his very famous and useful result that gives a canonical form of the fundamental matrix of \eqref{periodico_ordinario}:
\begin{theorem}{\textbf{(Floquet Theorem) (G. Floquet)} (\cite{Floquet})}\label{teo_floquet_clasico}

Let the ordinary homogeneous linear $\omega-$periodic differential system \eqref{periodico_ordinario}, where $A(t)$ is a continuous matrix.
Then, the fundamental matrix of system \eqref{periodico_ordinario} can be factorized in the Floquet form as 
    $X(t)=Q(t)\exp{(\Lambda t)},$
where $Q(t)$ is a $\omega-$periodic continuously differentiable matrix for $t\in\mathbb{R}$ and $\Lambda$ is a constant matrix.
\end{theorem}
The Floquet Theorem can be used to prove 
the following result stated by A.M. Lyapunov in his Ph.D. thesis (1892):
\begin{theorem}{\textbf{(Lyapunov reducibility theorem) (A.M. Lyapunov)} (\cite{ Lyapunov_tesis})}
Let the system \eqref{periodico_ordinario}, where $A(t)$ is a continuous matrix.
Then, system \eqref{periodico_ordinario} can be reduced to a system with constant coefficients by a linear non-singular continuous $\omega-$periodic Floquet-Lyapunov change of variables  $X=Q(t)Y$, transforming \eqref{periodico_ordinario} into the constant coefficients system $Y'(t)=\Lambda Y(t).$
\end{theorem}

The systems $X'(t)=A(t)X(t)$ and $Y'(t)=\Lambda Y(t)$ are \textit{Kinematically similar}. I.e., there exists a Lyapunov function $Q(t)$, satisfying $Q'(t)=A(t)Q(t)-Q(t)P$. In this case, $Q(t)$ is invertible, differentiable, and bounded (See \cite{daleckii_krein}).
The interested reader in periodic impulsive differential equations can see \cite{Bainov_periodico} and  \cite{Chicone, Cod-Lev, Feldman, Eastham_periodic} for further in Floquet theory for ordinary differential equations.\\

\noindent There is a remarkable quantity of literature about Floquet-Lyapunov theorems for another class of differential equations. We will present some relevant references concerning this work.\\
In \cite{Stokes} (1962), \textit{A. Stokes} gave an extension of the classical Floquet theorem class for the class of periodic functional differential equations
$x'_t(0)=f(x_t,t),$
where $x_t\in C,$ with $C$ is the space of continuous function defined from $[-h,0]$ to $\mathbb{R}^n,$ $h>0$, $A$ may be infinite, $x_t(\cdot)$ is defined as $x_t(s)=x(t+s), -h\leq s\leq 0$, $f:C\times \mathbb{R}\to \mathbb{R}^n,$ $f(\phi,t)$ linear in $\phi$, continuous and $\omega-$periodic satisfying $\Vert f(\phi,t)\Vert \leq L\vert \phi \vert,$ for some $L>0$ and $\forall (\phi,t),$ and $x_t'(0)$ denotes the right-had derivative of $x_t$ at $s=0$. I.e., $x_t'(0)=\lim_{r\to 0^+}\frac{1}{r}(x_{t+r}(0)-x_t(0)).$ 

In \cite{DACUNHA} (2011), \textit{Jeffrey J. DaCunha} and \textit{John M. Davis} studied
periodic linear systems on periodic time scales
$$x^{\Delta}(t)=A(t)x(t), \quad x(t_0)=x_0,$$
which include discrete, continuous, and mixed dynamical systems (hybrid dynamical systems). They gave a unified Floquet theorem that establishes a canonical Floquet decomposition on time scales in terms of the generalized exponential function and use these results to study homogeneous and nonhomogeneous periodic problems. 

In \cite{chandrika_floquet_DDE} (2023), \textit{J. Shaik, C. Prakash} and \textit{S. Tiwari} developed an approach to determining the stability of the following homogeneous linear $\omega$-periodic delay differential equation $x'(t)=a(t)x(t)+b(t)x(t-\tau)$,
where $a(t+\omega)=a(t), b(t+\omega)=b(t),$ and $x(t)=\eta(t), -\tau\leq t\leq 0,$ 
transforming the system into an approximate system of $\omega$-periodic ordinary differential equations using Galerkin approximations. Later, Floquet's theory is applied to the resultant ODEs. Since the original system is infinite-dimensional, they get an approximation by Floquet's normal solutions.\\

We emphasize that there is no literature on Floquet-Lyapunov theorems for DEPCA, IDEPCA, IDEPCA, DEPCAG, or IDEPCAG differential equations. Consequently, this seems to be the first work on this subject.  

%%%%%%%%%%%%%
\section{Aim of the work}
Inspired by A.M. Samoilenko and N.A. Perestyuk \cite{Samoilenko}, we will give a Floquet-Lyapunov type theorem for the class of nonautonomous homogeneous linear $\omega-$periodic \emph{IDEPCAG} 
\begin{equation}
\begin{tabular}{ll}
$x^{\prime }(t)=A(t)x(t)+B(t)x(\gamma (t)),$ & $t\neq t_{k},$ \\ 
$\Delta x|_{t=t_{k}}=C_{k}x(t_{k}^{-}),$ & $t=t_{k},$
\end{tabular}\label{periodico_idepcag_aim}
\end{equation}
with periodic conditions over all the coefficients involved. I.e., we will show that
\begin{itemize}
\item[(a)] The solutions of \eqref{periodico_idepcag_aim} can be represented in the \textbf{Floquet normal form} as
\begin{equation*}
X(t)=Q(t)\exp{(Pt)},\quad P=\dfrac{1}{\omega}Log\left(X(\omega)\right),\quad  t\in\mathbb{R},
\end{equation*}
where $P\in\mathbb{C}^{n\times n}$ is constant and the matrix function  $Q(t)\in \mathcal{PC}^{1}(\mathbb{R},\mathbb{C}^{n\times n})$ is  non-singular and $\omega-$periodic.
\item[(b)] System \eqref{periodico_idepcag_aim} can be reduced to the ordinary differential equation:
\begin{equation}
Y'(t)=PY(t),\label{DEPCAG_reducida_aim}
\end{equation}
by a $\omega-$periodic Floquet-Lyapunov transformation $X(t)=Q(t)Y(t).$ I.e the IDEPCAG \eqref{periodico_idepcag_aim} and  \eqref{DEPCAG_reducida_aim} are \textit{IDEPCAG-Kinematically similar} by the use of the Lyapunov function $Q(t)$, verifying the DEPCAG
$$Q'(t)=A(t)Q(t)-Q(t)P+B(t)Q(\gamma(t))e^{P(\gamma(t)-t)}.$$
\end{itemize}
%%%%%%%%%%%%%%%%

\section*{Why a Floquet theorem for IDEPCAG?}
Consider the following scalar IDEPCA
\begin{align}
x'(t)&=(A-1)x([t]),\qquad t\neq n, \nonumber \\
x(n)&=Cy(n^{-}),\qquad \qquad \,\,t=n, \quad n\in\mathbb{N} \label{idepca_parte_entera}.
\end{align}
where $A,C\in\mathbb{R}$ with $A,C\neq 1$ and $[t+1]=[t]+1,\,\forall t\in\mathbb{R}$. The equation \eqref{idepca_parte_entera} can be realized as an $1-$periodic system.\\
\noindent Let's solve \eqref{idepca_parte_entera}. If $t\in[n,n+1)$ for some $n\in\mathbb{Z}$, 
equation \eqref{idepca_parte_entera} can be written as
$x'(t)=(A-1)x(n).$\\
%\begin{equation*}
%    x'(t)=(A-1)x(n). \label{idepca_parte_entera_1}
%\end{equation*}
Without loss of generality, let $t_0=0$. Integrating on $[n,n+1)$ from $n$ to $t$, we get
\begin{equation}
    x(t)=x(n)(1+(A-1)(t-n)). \label{idepca_parte_entera_2}
\end{equation}
Next, assuming left-side continuity at $t=n+1$ and applying the impulse condition, we have
$x((n+1))=(AC)x(n).$ 
This is a \emph{finite-difference equation} whose solution is
\begin{equation}
    x(n)=(AC)^{n}x(0). \label{idepca_parte_entera_4}
\end{equation}
Finally, applying \eqref{idepca_parte_entera_4} in \eqref{idepca_parte_entera_2} we have found the solution of \eqref{idepca_parte_entera}
\begin{equation}
    x(t)=\left(AC\right)^{[t]}(1+(A-1)(t-[t]))x(0). \label{idepca_parte_entera_5}
\end{equation}
We can see that the nature of the dynamic is of mixed type. It depends on the discrete and the continuous parts of the system. The function $Q(t)=(1+(A-1)(t-[t]))$ is $1-$periodic and, from \eqref{idepca_parte_entera_5}, we can see the decomposition
\begin{equation*}
    x(t)=\exp{(Log(AC)[t])}\cdot (1+(A-1)(t-[t]))x_0,
\end{equation*}
suggests a Floquet normal form of the solution, where $Log(z),\,z\in\mathbb{C}-\{0\}$ is the principal complex logarithm. In this example, the presence of the impulse produces oscillations.
\begin{table}[h]
\centering
\begin{tabular}{|ll|}
\hline
\textbf{Behavior of solutions of \eqref{idepca_parte_entera}}  & \textbf{Condition}\\\hline
$x(t)$ is oscillatory and $x(t)\xrightarrow{t\to\infty} 0$  exponentially. & $-1<AC<0$ \\
$x(t)$ is nonoscillatory and $x(t)\xrightarrow{t\to\infty} 0$  exponentially. & $0<AC<1$ \\
$x(t)$ is nonoscillatory. & $AC\geq0$ \\
$x(t)$ is oscillatory. & $ AC<0$ \\
$x(t)$ is nontrivial $1-$periodic. & $AC=1$, with $A,C\neq 1,\,\,A,C>0$  \\
$x(t)$ is $2-$periodic and oscillatory. & $AC=-1$ with $A<0$ or $C<0$ \\ 
\hline
\end{tabular}
\end{table}
%%%%%%%%%%%%%%%
\begin{figure}[h!]
\centering
\includegraphics[scale=0.3]{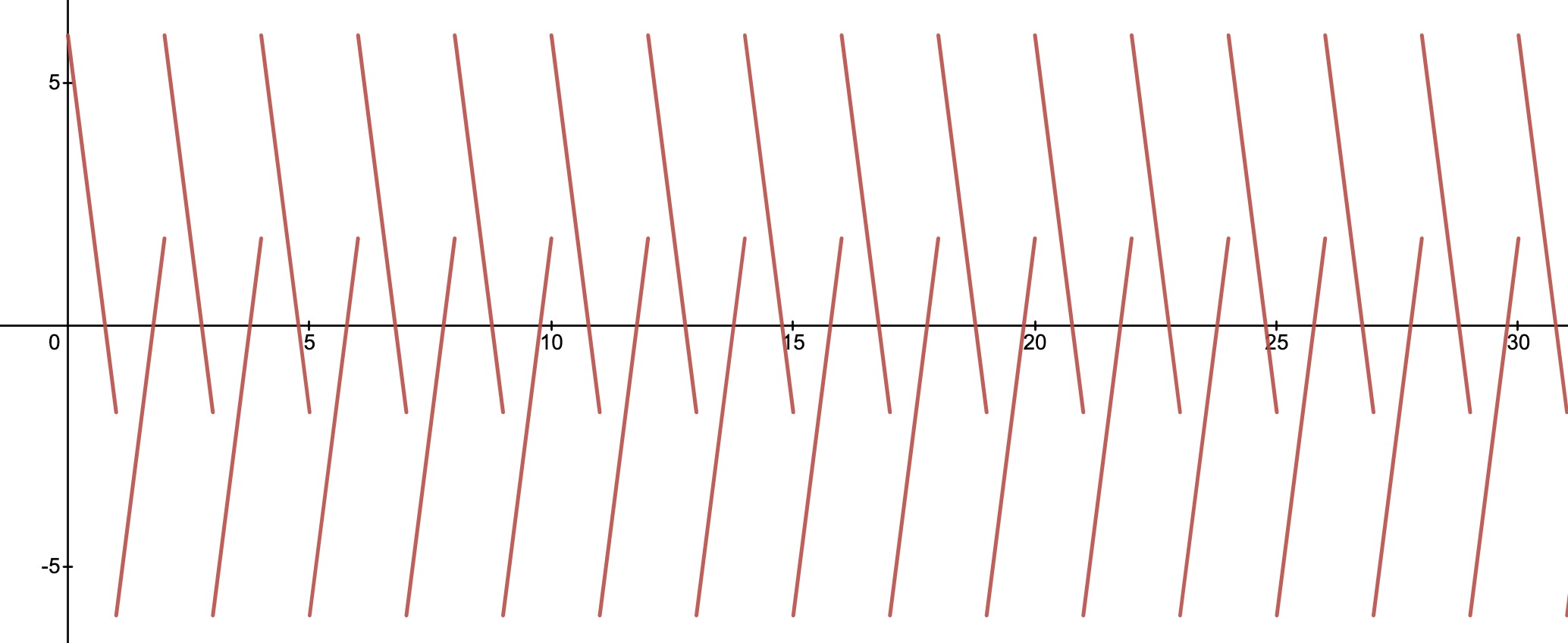}
\caption{Solution of \eqref{idepca_parte_entera} with
$A=-0.3 $, $C=10/3$, and $y_0=6$.}
 \label{ejemplos_1_y_2_periodicos}
\end{figure}

\section{Preliminaires of IDEPCAG}
Let $\mathcal{PC}(X, Y)$ be the set of all functions $r: X\to Y$ which are continuous for $t\neq t_k$ and continuous from the left with jump discontinuities at $t=t_k$. Similarly, let $\mathcal{PC}^{1}(X,Y)$ the set of functions $s:X\to Y$ such that $s'\in \mathcal{PC}(X,Y).$
\begin{definition}[DEPCAG solution]
A continuous function $x(t)$ is a solution of \eqref{depcag_eq} if:
\begin{itemize}
\item[(i)] $x'(t)$ exists at each point $t\in \mathbb{R}$ with the possible exception at the times $t_{k}$, $k\in \mathbb{Z}$, where the one side derivative exists.
\item[(ii)] $x(t)$ satisfies \eqref{depcag_eq} on the intervals of the form $(t_{k},t_{k+1})$, and it holds for the right
derivative of $x(t)$ at $t_{k}$.
\end{itemize}
\end{definition}
\begin{definition}[IDEPCAG solution]
A piecewise continuous function $z(t)$ is a solution of \eqref{idepcag_gral} if:
\begin{itemize}
\item[(i)] $z(t)$ is continuous on $I_k=[t_k,t_{k+1})$ with jump discontinuities at $t_k, \,k\in \mathbb{Z}$, 
where $z'(t)$ exists at each $t\in \mathbb{R}$ with the possible exception at the times $t_{k}$, where lateral derivatives exist (i.e. $z(t)\in \mathcal{PC}^{1}([t_k,t_{k+1}),R^n)$).
\item[(ii)] The ordinary differential equation
$$z^{\prime}(t)=f(t,z(t),z(\zeta_{k}))$$ holds on every interval $I_{k}$, where $\gamma(t)=\zeta_k$.
\item[(iii)]
For $t=t_{k}$, the impulsive condition 
$$\Delta z(t_{k})=z(t_{k})-z(t_{k}^{-})=J_{k}(z(t_{k}^{-}))$$
holds. I.e., $z(t_{k})=z(t_{k}^{-})+J_{k}(z(t_{k}^{-}))$, where $z(t_{k}^{-})$ denotes the left-hand limit
of the function $y$ at $t_k$.
\end{itemize}
\end{definition}
%%%%%%
\subsection{\textbf{Solving the nonautonomous homogeneous linear IDEPCAG}}
In this section, we will present the nonautonomous homogeneous linear IDEPCAG 
\begin{equation}
\begin{tabular}{ll}
$x^{\prime }(t)=A(t)x(t)+B(t)x(\gamma (t)),$ & $t\neq t_{k}$ \\ 
$\Delta x|_{t=t_{k}}=C_k x(t_{k}^{-}),$ & $t=t_{k}$
\end{tabular}
\label{SISTEMA_IDEPCAG_ORIGINAL_}
\end{equation}
where $x\in\mathbb{C}^n,t\in\mathbb{R},$ $A(t),C(t)$ are real-valued continuous locally integrable $n\times n$ matrix functions, $(C_k)_{k\in\mathbb{N}}$ is a real $n\times n$ matrix sequence such that $\det(I+C_k)\neq 0$ $\forall k\in\mathbb{N},$ where $I$ is the $n\times n$ identity matrix and $\gamma (t)$ is a generalized piecewise constant argument.\\

\noindent During the rest of the work, we will assume $\gamma(\tau):=\tau$ if $t_{k(\tau)}\leq \gamma(\tau)< \tau <t_{k(\tau)+1},$ where $k(\tau)$ is the only $k\in\mathbb{Z}$ such that $t_{k(\tau)}\leq \tau \leq t_{k(\tau)+1}.$\\

\noindent Let $z(t)=\Phi(t,\tau)z(\tau)$
the  solution of the ordinary differential equation
\begin{eqnarray*}
z^{\prime}(t)=A(t)z(t),\quad z_0=z(\tau), \qquad t,\tau\in [\tau,\infty),
\end{eqnarray*}
where $\Phi(t,s)=\Phi(t)\Phi^{-1}(s),\,\,\Phi(t,u)\Phi(u,s)=\Phi(t,s).$\\
For the sake of simplicity, we will consider the normalized fundamental matrix $\Phi(0)=I.$ All our results can be rewritten considering an arbitrary value of $\Phi(0).$\\

We will assume the following hypothesis:
\begin{enumerate}
    \item[\textbf{(H)}] Let
$$\displaystyle{\sigma_{k}^{+}(A) =exp\left(\int_{t_{k}}^{\zeta_{k}}\left| A(u)\right| du\right),}\qquad  
\displaystyle{\sigma_{k}^{-}(A)
=exp\left(\int_{\zeta_{k}}^{t_{k+1}}\left| A(u)\right| du\right),}$$ 
$$\sigma_{k}(A)=\sigma_{k}^{+}(A)\sigma_{k}^{-}(A),\qquad
\nu_{k}^{\pm }(B)=\sigma_{k}^{\pm }(A)\ln \sigma_{k}^{\pm }(B),$$

and assume that $$\sigma(A)=\displaystyle{\sup_{k\in \mathbb{Z}}\sigma_{k}(A)<\infty},\qquad \nu^{\pm}(B)=\sup_{k\in \mathbb{Z}}\nu_{k}^{\pm}(B)<\infty,$$
where
\begin{equation}
\nu_k^{+}(B)<\nu^{+}(B)<1,\quad \nu_k^{-}(B)<\nu^{-}(B)<1. \label{Cond_invertibilidad}
\end{equation}
\end{enumerate}
\noindent Consider the following definitions
\begin{equation}
J(t,\tau)=I+\int_{\tau}^{t}\Phi(\tau,s)B(s)ds,\quad
E(t,\tau)=\Phi(t,\tau) J(t,\tau), \label{MATRIZ J}
\end{equation}
where $I$ is the $n\times n$ identity matrix and $|\cdot|$ is some matricial norm.
\begin{remark}
As a consequence of $\textbf{(H)}$, it is important to notice the following facts: 
\begin{itemize}
\item[(i)] Due to condition \eqref{Cond_invertibilidad},
$J^{-1}(t_k,\zeta_{k})$ and $J^{-1}(t_{k+1},\zeta_{k})$ are well defined $\forall k\in \mathbb{Z},$ and 
\begin{eqnarray*}
\left|J^{-1}(t_k,\zeta_k)\right| \leq \sum\limits_{k=0}^{\infty } \left[ \nu^{+} (b)\right] ^{k}  =\frac{1}{1-\nu^{+} (b)},\quad 
\left| J(t_{k+1},\zeta_{k})\right| \leq 1+\nu^{-} (b), \label{COTA_J}\\
\left|J^{-1}(t_{k+1},\zeta_k)\right| \leq \sum\limits_{k=0}^{\infty } \left[ \nu^{-} (b)\right] ^{k}  =\frac{1}{1-\nu^{-} (b)},\quad 
\left| J(t_{k},\zeta_{k})\right| \leq 1+\nu^{+} (b).\label{COTA_J2}
\end{eqnarray*}
Also, if we set $t_0=\tau$, we are considering that $J^{-1}(\tau,\gamma(\tau))$ exists.

\item[(ii)] On the other hand, if we want a non-zero solution of a linear IDEPCAG, we need $J(t,\zeta_k)\neq 0, \forall k\in\mathbb{Z}$ and $\forall t\in [\tau,\infty)$ (See Remark \ref{remark_ceros} and \cite{Kuo_pinto_2013}).
\end{itemize}
\end{remark}
%%%%%%%%%%%%%%%
In the rest of this work, we will also assume the following notation:
\begin{equation*}
\prod_{j=1}^{n}A_j=
    \begin{cases}
        A_n\cdot A_{n-1} \cdots A_1, & \text{if } n\geq 1,\\
        \qquad \quad I & \text{if } n<1.
    \end{cases}
    \quad \text{ and } \quad
\sum_{j=1}^{n}A_j=
    \begin{cases}
        A_1+\ldots +A_{n}, & \text{if } n\geq 1,\\
        \qquad \quad 0 & \text{if } n<1.
    \end{cases}
\end{equation*}
\begin{remark}
Also, for writing and space convenience, we will denote the right-side matricial product of $A$ and  $B^{-1}$ as
$A\cdot B^{-1}=\dfrac{A}{B}$.    
\end{remark}

\subsection{The fundamental solution of the homogeneous linear IDEPCAG}
The following results can be found in  \cite{Torres1} and \cite{TORRES_var}. They are the IDEPCAG extension of \cite{P2011} (the case with $C_k=0\,,\forall k\in\mathbb{Z})$:
\begin{theorem}\label{TEO_FORMULA_var_PAram}\cite{TORRES_var}
Let the following linear IDEPCAG system
\begin{equation}
\begin{tabular}{ll}
$X^{\prime }(t)=A(t)X(t)+B(t)X(\gamma (t)),$ & $t\neq t_{k}$ \\ 
$X(t_{k})=\left(I+C_k\right) X(t_{k}^{-}),$ & $t=t_{k}$\\
$X_0=X(\tau).$
\end{tabular}
\label{sistema_w}
\end{equation}
If \textbf{(H)} holds, then the unique solution of \eqref{sistema_w} is
 \begin{equation}
X(t)=W(t,\tau)z(\tau),  \quad t\in [\tau,\infty),\label{SOLUCION_FINAL_SISTEMA_IDEPCAG_LINEAL}
%\label{DEFINICION_GENERAL_MATRIZ_W}
\end{equation}
where $W(t,\tau)$ is given by 
\begin{equation}
W(t,\tau)=W(t,t_{k(t)})\left(\prod_{r=k(\tau)+2}^{k(t)}\left(I+C_{r}\right) W(t_{r},t_{r-1})\right) 
\left(I+C_{k\left(\tau \right) +1}\right) W(t_{k(\tau)+1},\tau )
\label{MATRIZ_FUNDAMENTAL_IDEPCAG}
\end{equation}
for $t\in I_{k(t)}, \tau \in I_{k(\tau) },$
and $W(t,s)$ is defined as 
\begin{equation*}
W(t,s)=\dfrac{E(t,\gamma(s))}{E(s ,\gamma(s))},\qquad \text{if }t,s
\in I_k=[t_{k},t_{k+1}]. \label{MATRIZ_w}
\end{equation*}
Also, the discrete solution of \eqref{sistema_w} is given by 
\begin{equation}
X(t_{k(t)})=\left(\prod_{r=k(\tau)+2}^{k(t)}\left(I+C_{r}\right) W(t_{r},t_{r-1})\right) \left(I+C_{k(\tau)+1}\right) W(t_{k(\tau)+1},\tau)X(\tau).
\label{SOLUCION_DISCRETA_SISTEMA_Z}
\end{equation}
The expression \eqref{MATRIZ_FUNDAMENTAL_IDEPCAG}
is called the Cauchy matrix of \eqref{sistema_w}.
\end{theorem}
\begin{proof}
    Let $t,\tau \in I_k=[t_{k},t_{k+1})$ for some $k\in \mathbb{Z}.$
In this interval, we are in the presence of the ordinary system
$$X^{\prime }(t)=A(t)X(t)+B(t)X(\zeta_k).$$
So, the unique solution can be written as
\begin{equation}
X(t)=\Phi(t,\tau)X(\tau )+\int_{\tau}^{t}\Phi(t,s)B(s)X(\zeta_{k})ds.
\label{variacion_parametros_general}
\end{equation}
Keeping in mind $I_k^{+}$, evaluating the last expression at  $t=\zeta_{k}$ we have
\begin{equation*}
X(\zeta_{k})=\Phi(\zeta_{k},\tau )X(\tau)+\int_{\tau}^{\zeta_{k}}\Phi(\zeta_{k},s)B(s)X(\zeta_{k})ds.
%\label{variacion_parametros_gamma_general}
\end{equation*}
Hence, we get
\begin{eqnarray*}
\left(I+\int_{\zeta_{k}}^{\tau}\Phi(\zeta_{k},s)B(s)ds\right)X(\zeta_{k}) &=&\Phi(\zeta_{k},\tau )X(\tau),
\end{eqnarray*}
i.e
\begin{equation*}
X(\zeta_{k})=J^{-1}(\tau ,\zeta_{k})\Phi(\zeta_{k},\tau)X(\tau).
\label{inicial_gamma_general}
\end{equation*}
Then, by the definition of $E(t,\tau)=\Phi(t,\tau)J(t,\tau)$, we have
\begin{equation}
X(\zeta_{k})=E^{-1}(\tau ,\zeta_{k})X(\tau).
\label{condicion_inicial_con_gama_y_E}
\end{equation}
Now, from \eqref{variacion_parametros_general} working on $I_k^{-}$, considering $\tau =\zeta_{k}$, we have
\begin{eqnarray*}
X(t)&=&\Phi(t,\zeta_{k})X(\zeta_{k})+\int_{\zeta_{k}}^{t}\Phi(t,s)B(s)X(\zeta_{k})ds \\
&=&\Phi(t,\zeta_{k})\left( I+\int_{\zeta_{k}}^{t}\Phi(\zeta_{k},s)B(s)ds\right) X(\zeta_{k}),
\end{eqnarray*}
i.e.,
\begin{equation}
X(t)=E(t,\zeta_{k})X(\zeta_{k}).
\label{lineal_homogenea_parte_retardo_sin_cond_inicial}
\end{equation}
So, by \eqref{condicion_inicial_con_gama_y_E}, we can rewrite  \eqref{lineal_homogenea_parte_retardo_sin_cond_inicial} as
\begin{equation}
X(t)=\dfrac{E(t,\zeta_{k})}{E(\tau ,\zeta_{k})}X(\tau).
\label{condicion_inicial_sin_w}
\end{equation}

\noindent Then, setting  
\begin{equation*}
W(t,s)=\dfrac{E(t,\gamma(s))}{E(s ,\gamma(s))},\qquad \text{if }t,s
\in I_k=[t_{k},t_{k+1}],  \label{MATRIZ_w_}
\end{equation*}
we have the solution for \eqref{sistema_w} for $t\in I_{k}=[t_{k},t_{k+1}),$
\begin{equation}
X(t)=W(t,\tau)X(\tau). \label{DEFINICION_GENERAL_MATRIZ_W_}
\end{equation}

\noindent Next, if we consider $\tau=t_k, $ and, assuming left side continuity of \eqref{SOLUCION_FINAL_SISTEMA_IDEPCAG_LINEAL} at $t=t_{k+1}$, we have
\begin{equation*}
X(t_{k+1}^{-})=W(t_{k+1},t_k)X(t_k)
\end{equation*}
Then, applying the impulsive condition defined in \eqref{sistema_w} to the last equation, we get
\begin{eqnarray}
X(t_{k+1}) &=&\left(I+C_{k+1}\right)W(t_{k+1},t_k)X(t_k).\label{cond_discreta_osc_kuo}
\end{eqnarray}

\noindent The last expression defines a finite-difference equation whose solution is \eqref{SOLUCION_DISCRETA_SISTEMA_Z}. 
Now, by \eqref{DEFINICION_GENERAL_MATRIZ_W_} and the impulsive condition defined in \eqref{sistema_w}, we have
$$X(t_{k(\tau)+1})=(I+C_{k(\tau)+1})W(t_{k(\tau)+1},\tau)X(\tau).$$
Hence, considering $\tau=t_k$ in \eqref{SOLUCION_FINAL_SISTEMA_IDEPCAG_LINEAL} and applying \eqref{SOLUCION_DISCRETA_SISTEMA_Z}, 
we get \eqref{SOLUCION_FINAL_SISTEMA_IDEPCAG_LINEAL}.
In this way, we have solved \eqref{sistema_w} on $[\tau,t).$\\
We used the decomposition of $I_k=I_k^+\cup I_k^-$ to define $W$. In fact, we can rewrite \eqref{MATRIZ_FUNDAMENTAL_IDEPCAG} in terms of the advanced and delayed parts using \eqref{MATRIZ_w}:
\begin{eqnarray}
W(t,\tau)&=&\dfrac{E(t,\zeta_{k(t)})}{E(t_{k(t)},\zeta_{k(t)})}\left(\prod_{r=k(\tau)+2}^{k(t)}\left(I+C_{r}\right) \dfrac{E(t_{r},\zeta_{r-1})}{E(t_{r-1},\zeta_{r-1})}\right)\label{MATRIZ_FUNDAMENTAL_IDEPCAG_avance_retardo_1}\\ 
&&\cdot \left(I+C_{k\left(\tau \right)+1}\right) \dfrac{E(t_{k(\tau)+1},\gamma(\tau))}{E(\tau,\gamma(\tau))},\qquad \zeta_r=\gamma(t_r),\notag
\end{eqnarray}
for $t\in I_{k(t)}$ and $\tau \in I_{k(\tau) }.$
\end{proof}
\begin{remark}\label{remark_ceros}
\begin{itemize}
\item[]
\item[(i)]   Considering $B(t)=0$, we recover the classical fundamental solution of the impulsive linear differential equation (see \cite{Samoilenko}).
\item[(ii)] If $C_k=0, \forall k\in\mathbb{Z}$, we recover the DEPCAG case studied by M. Pinto in \cite{P2011}. 
\end{itemize}
\end{remark}

\section{The Floquet theory for IDEPCAG}
Let the $\omega-$periodic homogeneous linear IDEPCAG 
\begin{equation}
\begin{tabular}{ll}
$X^{\prime }(t)=A(t)X(t)+B(t)X(\gamma (t)),$ & $t\neq t_{k}$ \\ 
$\Delta X|_{t=t_{k}}=C_{k}X(t_{k}^{-}),$ & $t=t_{k}$%
\end{tabular}
\label{sistema_idepcag_periodico_homogeneo}
\end{equation}
where $A(t), B(t)$ are continuous $n\times n$ real-valued locally integrable matrix functions (piecewise continuous with jump discontinuities at $t=t_k$), and there exists a natural number $p$ such that $\det(I+C_k)\neq 0,\forall k=1,2,\ldots, p$ and
\begin{eqnarray}
    &&A(t+\omega)=A(t),\quad B(t+\omega)=B(t),\quad \forall t\in[0,\infty),\nonumber\\
    &&C_{k+p}=C_k,\quad \forall k\in\mathbb{Z},\label{Periodicidad_coeficientes}
\end{eqnarray}
$t_0=\tau<t_1<\ldots<t_p\leq \tau+\omega,$  and $\gamma$ is a piecewise constant argument of generalized type such that
$\gamma(t)=\zeta_{k}$ if $t\in[t_k,t_{k+1})$ with $t_k\leq \zeta_k\leq t_{k+1},$ with the so-called $(\omega,p)-$property
\begin{eqnarray}
    &&t_{k+p}=t_k+\omega,\qquad \zeta_{k+p}=\zeta_k+\omega, \quad  \forall k\in\mathbb{Z}.\label{w-p_propiedad} 
\end{eqnarray}

This section will provide an IDEPCAG version of the Floquet Theorem.
\subsection{\textbf{Auxiliary results}}
In the following, we will assume the classical Floquet Theorem for the solutions of the $\omega-$periodic ordinary system
\begin{eqnarray}
Z'(t)=A(t)Z(t),\label{homogeneo_periodico_ordinaro}\\
A(t+\omega)=A(t),\nonumber
\end{eqnarray}
with $\Phi(\tau)=I$. I.e., $\Phi(t+\omega)=\Phi(t)\Phi(\omega), \,\forall t\in\mathbb{R}.$ 
\begin{lemma}\label{periodicidad_j_phi_e}
Let the matrices $J(t,s),\,\Phi(t,s)$ and $E(t,s)$ as they were defined on $\textbf{(H)}$. Then, the following properties hold:
\begin{equation}
\Phi(t+\omega,s+\omega)=\Phi(t,s),\,\,
J(t+\omega,s+\omega)=J(t,s),\,\,
E(t+\omega,s+\omega)=E(t,s),\,\,\forall t,s\in\mathbb{R}.\label{biperiodicidad_operadores}
\end{equation}
\end{lemma}
\begin{proof}
Because the classical Floquet Theorem applied \eqref{homogeneo_periodico_ordinaro}, we have $\Phi(t+\omega)=\Phi(t)\Phi(\omega).$ Then
\begin{eqnarray*}
\Phi(t+\omega,s+\omega)&=&\Phi(t+\omega)\Phi^{-1}(s+\omega)\\
&=&\Phi(t)\Phi(\omega)\Phi^{-1}(\omega)\Phi^{-1}(s)\\
&=&\Phi(t)\Phi^{-1}(s)\\
&=&\Phi(t,s).
\end{eqnarray*}
Next, in order to prove the biperiodicity of $J(t,s)$, using the $\omega-$periodicity of $B(t)$, the $\omega-$biperiodicity of $\phi(t,s)$ and the change of variables $z=u-\omega$, we see that
\begin{eqnarray*}
    J(t+\omega, s+\omega)&=&\displaystyle{I+\int_{s+\omega}^{t+\omega}\Phi(s+\omega,u)B(u)du}\\
    &=&\displaystyle{I+\int_{s+\omega}^{t+\omega}\Phi(s+2\omega,u+\omega)B(u)du}\\
    &=&\displaystyle{I+\int_{s+\omega}^{t+\omega}\Phi(s+2\omega,u+\omega)B(u-\omega)du}\\
    &=&\displaystyle{I+\int_{s}^{t}\Phi(s+2\omega,z+2\omega)B(z)dz}\\
    &=&I+\displaystyle{\int_{s}^{t}}\Phi(s,z)B(z)dz\\
    &=&J(t,s).
\end{eqnarray*}
\noindent Hence, as $E(t,s)=\Phi(t,s)J(t,s)$,  we also conclude that $E(t+\omega,s+\omega)=E(t,s).$
\end{proof}
As a corollary, using \eqref{MATRIZ_FUNDAMENTAL_IDEPCAG_avance_retardo_1}, it is easy to prove the following result:
\begin{corollary}
Let Lemma \ref{periodicidad_j_phi_e} holds. Then, the so-called Transition matrix (matriciant or Cauchy matrix) associated with \eqref{sistema_idepcag_periodico_homogeneo} satisfies $W(t+\omega,s+\omega)=W(t,s),\,\, \forall t,s\in\mathbb{R}.$
\end{corollary}

\subsection{\textbf{The Monodromy operator}}
Some of the following are basic results; nevertheless, we will present them for a better understanding and completeness. They can be found at \cite{Bainov_periodico,Brown2013}:
\begin{lemma}
    If $X(t)$ is a fundamental solution of \eqref{sistema_idepcag_periodico_homogeneo}, then $X(t+\omega)$ also is a fundamental matrix of \eqref{sistema_idepcag_periodico_homogeneo}.
\end{lemma}
\begin{proof}
Let $Y(t)=X(t+\omega).$ Then, for $t\neq t_k$, we have
    \begin{eqnarray*}
Y'(t)&=&A(t+\omega)Y(t)+B(t+\omega)Y(\gamma(t))\\
&=&A(t)Y(t)+B(t)Y(\gamma(t)).
    \end{eqnarray*}
Finally, for $t=t_k$ and setting $Y(t_k)=X(t_k+\omega)=X(t_{k+p}),$ we have
\begin{eqnarray*}
    \Delta Y(t_k)=\Delta X(t_{k+p})=C_{k+p}X(t^-_{k+p})=C_kY(t_k).
\end{eqnarray*}
\end{proof}
\noindent 
%Without loss of generality, let $t_0=0$ and 
Let $\zeta_j=\gamma(t_j),\forall j\in\mathbb{Z}$ and define $\zeta_0=\zeta_{k(\tau)}\coloneqq \gamma(\tau),\quad  t_0\coloneqq t_{k(\tau)}=\tau.$ Since the $(\omega,p)$-property \eqref{w-p_propiedad} and Lemma \ref{periodicidad_j_phi_e}, we have $t_{k(\tau+\omega)}=\tau+\omega,\,\, \zeta_{k(\tau+\omega)}=\gamma(\tau)+\omega,$
and $$\dfrac{E(\tau+\omega,\zeta_{k(\tau+\omega)})}{E(t_{k(\tau+\omega)},\zeta_{k(\tau+\omega)})}
=\dfrac{E(\tau+\omega,\zeta_{k(\tau)}+\omega)}{E(t_{k(\tau)}+\omega,\zeta_{k(\tau)}+\omega)}
%=\dfrac{E(\tau+\omega,\gamma(\tau)+\omega)}{E(\tau+\omega,\gamma(\tau)+\omega)}
=I.$$
Therefore, if we consider $X(\tau)=I$ and evaluating at $t=\tau+\omega$ in \eqref{MATRIZ_FUNDAMENTAL_IDEPCAG_avance_retardo_1}, we have
\begin{eqnarray*}
X(\tau+\omega)&=&\dfrac{E(\tau+\omega,\zeta_{k(\tau+\omega)})}{E(t_{k(\tau+\omega)},\zeta_{k(\tau+\omega})}\left(\prod_{r=1}^{p}\left(I+C_{r}\right) \dfrac{E(t_{r},\zeta_{r-1})}{E(t_{r-1},\zeta_{r-1})}\right)=X(\omega).
\end{eqnarray*}
Hence, we can define
\begin{equation}
X(\omega)=\prod_{r=1}^{p}\left(I+C_{r}\right) \dfrac{E(t_{r},\zeta_{r-1})}{E(t_{r-1},\zeta_{r-1})}\\
\label{MATRIZ_monodromia_avance_retardo}
\end{equation}
as the so-called \textbf{monodromy operator} or \textbf{monodromy matrix} of \eqref{sistema_w}. 
Notice that we have shown  $X(\tau+\omega)=X(\tau)X(\omega),$ where $X(\tau)=I.$\\

Without loss of generality, in the rest of the work, we will consider $t_0=\tau=0.$
\begin{theorem}{(\textbf{Floquet factorization theorem})}\label{Teo_forma_normal}

\noindent Let Lemma \ref{periodicidad_j_phi_e} holds. Then, the fundamental solution of \eqref{sistema_idepcag_periodico_homogeneo} $X(t)$ with $X(0)=I$ can be written in the Floquet normal form as 
\begin{equation}
X(t+\omega)=X(t)X(\omega).\label{factorizacion_floquet}
\end{equation}
\end{theorem}    
\begin{proof}
We will compute $X(t+\omega)$ directly. Evaluating \eqref{MATRIZ_FUNDAMENTAL_IDEPCAG_avance_retardo_1} at $t+\omega$, we get  
\begin{eqnarray}
X(t+\omega)&=&\dfrac{E(t+\omega,\zeta_{k(t)+p})}{E(t_{k(t)+p},\zeta_{k(t)+p})}\left(\prod_{r=1}^{k(t)+p}\left(I+C_{r}\right) \dfrac{E(t_{r},\zeta_{r-1})}{E(t_{r-1},\zeta_{r-1})}\right)\nonumber\\
&=&\dfrac{E(t+\omega,\zeta_{k(t)+p})}{E(t_{k(t)+p},\zeta_{k(t)+p})}\left(\prod_{r=p+1}^{k(t)+p}\left(I+C_{r}\right) \dfrac{E(t_{r},\zeta_{r-1})}{E(t_{r-1},\zeta_{r-1})}\right)\left(\prod_{r=1}^{p}\left(I+C_{r}\right) \dfrac{E(t_{r},\zeta_{r-1})}{E(t_{r-1},\zeta_{r-1})}\right)\nonumber\\
&=&\dfrac{E(t+\omega,\zeta_{k(t)+p})}{E(t_{k(t)+p},\zeta_{k(t)+p})}\left(\prod_{r=1}^{k(t)}\left(I+C_{r+p}\right) \dfrac{E(t_{r+p},\zeta_{r-1+p})}{E(t_{r-1+p},\zeta_{r-1+p})}\right)X(\omega)\nonumber\\
&=&\dfrac{E(t+\omega,\zeta_{k(t)}+\omega)}{E(t_{k(t)}+\omega,\zeta_{k(t)}+\omega)}\left(\prod_{r=1}^{k(t)}\left(I+C_{r}\right) \dfrac{E(t_{r}+\omega,\zeta_{r-1}+\omega)}{E(t_{r-1}+\omega,\zeta_{r-1}+\omega)}\right)X(\omega)\nonumber\\
&=&\dfrac{E(t,\zeta_{k(t)})}{E(t_{k(t)},\zeta_{k(t)})}\left(\prod_{r=1}^{k(t)}\left(I+C_{r}\right) \dfrac{E(t_{r},\zeta_{r-1})}{E(t_{r-1},\zeta_{r-1})}\right)X(\omega)\nonumber\\
&=&X(t)X(\omega).\nonumber
\end{eqnarray}
\end{proof}
As a consequence of the last Theorem, we have a necessary and sufficient condition for the existence of an $\omega-$ periodic solution for the IDEPCAG \eqref{sistema_idepcag_periodico_homogeneo}:
\begin{corollary}{\textbf{(Criterion for existence of periodic solutions for IDEPCAG \eqref{sistema_idepcag_periodico_homogeneo})}}\label{Existencia_soluciones_periodicas}

Let the fundamental solution of \eqref{sistema_idepcag_periodico_homogeneo} $X(t)$ with $X(0)=I$ and  Lemma \ref{periodicidad_j_phi_e} holds. Then, \eqref{sistema_idepcag_periodico_homogeneo} has an $\omega-$periodic solution if and only if $X(\omega)=I.$ I.e.,
     \begin{eqnarray}
\prod_{r=1}^{p}\left(\left(I+C_{r}\right) \dfrac{E(t_{r},\zeta_{r-1})}{E(t_{r-1},\zeta_{r-1})}\right)=I.
\end{eqnarray}
\end{corollary}
\begin{proof}
Let Theorem \ref{Teo_forma_normal} holds. 
\begin{itemize}
        %\item [] 
\item[$(\Leftarrow)$] If $X(\omega)=I,$ we have
\begin{eqnarray*}
X(t+\omega)&=&X(t)X(\omega)\\
&=&X(t).
\end{eqnarray*}
\item[$(\Rightarrow$)] If $X(t+\omega)=X(t),$ then, evaluating at $t=0$, we have $X(\omega)=X(0)=I.$
\end{itemize}
\end{proof}
\begin{corollary}\label{Existencia_soluciones_KW_periodicas}
     Let the $N\in\mathbb{N}$. Let the fundamental solution of \eqref{sistema_idepcag_periodico_homogeneo} $X(t)$ with $X(0)=I$ and  Lemma \ref{periodicidad_j_phi_e} holds. Then, \eqref{sistema_idepcag_periodico_homogeneo} has an $N\omega-$periodic solution if and only if $X^{N}(\omega)=I,$ where $I$ is the identity matrix. I.e.,
     \begin{eqnarray*}
\left(\prod_{r=1}^{p}\left(\left(I+C_{r}\right) \dfrac{E(t_{r},\zeta_{r-1})}{E(t_{r-1},\zeta_{r-1})}\right)\right)^N=I.
\end{eqnarray*}
\end{corollary}

\begin{remark}
\begin{enumerate}
    \item[]
    \item Because of \eqref{biperiodicidad_operadores}, the fundamental matrix of a homogeneous linear $\omega-$periodic IDEPCAG system already has the Floquet factorization form given by \eqref{factorizacion_floquet}. This is a remarkable and expected fact.
    \item Corollary \ref{Existencia_soluciones_periodicas} is an extension of the condition given by K-S. Chiu and M. Pinto in \cite{Kuo_pinto_2013} for the existence of $\omega-$periodic solutions of homogeneous linear DEPCAG case. The authors considered $t_p=\omega$ and $C_j=0,\,\forall j\in\mathbb{Z}.$.
\end{enumerate}
\end{remark}

\subsection{\textbf{The Logarithm of the monodromy operator}}
As indicated before, we will consider
$Log(z)$ as the complex principal logarithm with $$Log(z)=\ln(|z|) +i\arg(z),\,\,\,-\pi<\arg(z)\leq \pi \text{ and }z\neq 0.$$
In this section, we will give some conditions for the existence of a logarithm of a matrix.
%%%%%%%%%%%%%%%%%%
\subsection{\textbf{Floquet Multipliers, Floquet exponents and Lyapunov exponents}}

\subsubsection{\textbf{Floquet multipliers}}
\begin{definition}
    The eigenvalues $\rho_1,\rho_2,\ldots, \rho_n$ (counting multiplicities) of the Monodromy matrix $X(\omega)$ are the so-called \textbf{Floquet multipliers of $X(\omega)$.}
\end{definition} 
We know that the Floquet multipliers are non-zero since $X(t+\omega)$ and $X(t)$ are fundamental matrices of \eqref{sistema_idepcag_periodico_homogeneo}, and therefore, non-singular. In fact, 
\begin{equation}
\det(X(\omega))=\dfrac{\det(X(t+\omega))}{\det(X(t))}=\displaystyle{\prod_{i=1}^n \rho_i}\neq 0. \label{monodromia_determinante}   
\end{equation}
As $\rho_j\neq 0, \,\forall j\in\{1,2,\ldots, n\},$ we can write the Floquet multipliers as $$\rho_j=\exp{(\lambda_j)},\,\,\lambda_j\in\mathbb{C}.$$
An amazing fact is that the dynamics of the $\omega-$periodic system \eqref{sistema_idepcag_periodico_homogeneo} is governed by the spectral properties of $X(\omega)$. The Floquet multipliers  will play a crucial role in that purpose:
\begin{theorem}\label{floquet_multiplier}
    Let Theorem \ref{Teo_forma_normal} holds and consider the Monodromy matrix $X(\omega)$ of the $\omega-$periodic system \eqref{sistema_idepcag_periodico_homogeneo}. Then, a Floquet multiplier $\rho_j=\exp{(\lambda_j)}$ with $\lambda_j\in\mathbb{C}$ is an eigenvalue of $X(\omega)$ if and only if there is a non-trivial solution $x_j:\mathbb{R}\to\mathbb{C}$ such that $$x_j(t+\omega)=\rho_j x_j(t),\quad  j\in\{1,2,\ldots,n\}, \,t\in\mathbb{R}.$$
\end{theorem}
\begin{proof}
    Let $v_j\in\mathbb{C}^n-\{0\}$ be an eigenvector of $X(\omega)$ for the eigenvalue $\rho_j=\exp{(\lambda_j)},$ and set $$x_j(t)\coloneqq X(t)v_j,$$ where $X(t)$ is the fundamental matrix of \eqref{sistema_idepcag_periodico_homogeneo} with $X(0)=I.$ Then, $x_j(t)$ is a solution of \eqref{sistema_idepcag_periodico_homogeneo} and
    \begin{eqnarray*}
        x_j(t+\omega)&=&X(t+\omega)v_j\\
        &=&X(t)X(\omega)v_j\\
        &=&\rho_j X(t)v_j\\
        &=&\rho_j x_j(t).
    \end{eqnarray*}
    Conversely, if $x_j(t):\mathbb{R}\to\mathbb{C}$ is a nontrivial solution satisfying $x_j(t+\omega)=\rho_j x_j(t),$ we can consider $x_j(0)\neq 0$. Then, we see that 
    \begin{eqnarray*}  x_j(\omega)=X(\omega)x_j(0)=\rho_j x_j(0).
    \end{eqnarray*}
I.e., $x_j(0)$ is an eigenvector of $X(\omega)$  with associated eigenvalue $\rho_j.$
\end{proof}
It is important to remark that if $Y(t)$ is any other fundamental matrix for \eqref{sistema_idepcag_periodico_homogeneo}, then
    $$X(t)=Y(t)G,$$ for some non-singular matrix $G$. So, we can see that:
\begin{eqnarray*}
    Y(t+\omega)G&=&X(t+\omega)\\
    &=&X(t)X(\omega)\\
    &=&Y(t)G X(\omega).
\end{eqnarray*}
    I.e.,  $Y(t+\omega)=Y(t)GX(\omega)G^{-1}.$
    Hence, by the last equation, every fundamental matrix $Y(t)$ determines a matrix $GX(\omega)G^{-1}$. Since, as the spectrum of $X(\omega)$ is invariant under similarity, all the fundamental matrices have the same  Floquet multipliers.\\

As a corollary of Theorem \ref{floquet_multiplier}, we have the following result concerning the asymptotic behavior of the solutions of \eqref{sistema_idepcag_periodico_homogeneo}:
\begin{corollary}{\textbf{(Asymptotic behavior of the solutions of a $\omega-$periodic linear IDEPCAG by Floquet multipliers})}\label{teo_asympt_behavior_floquet_multipliers}

The solutions of \eqref{sistema_idepcag_periodico_homogeneo} converges exponentially to zero if $|\rho_j|<1$, they will be  $\omega-$periodic (or $2\omega$-periodic) if $|\rho_j|=1$ and they will be unbounded if $|\rho_j|>1.$ In other words, if the Floquet multipliers lie in the unit circle, solutions of \eqref{sistema_idepcag_periodico_homogeneo} will be bounded. Otherwise, they will be unbounded.  
\end{corollary}

\subsubsection{\textbf{Floquet exponents}}
\begin{definition}
    Let $\rho_j=\exp{(\lambda_j)}, \,j\in\{1,2,\ldots,n\}$ a Floquet multiplier of $X(\omega)$. We will call to the number $\dfrac{1}{\omega}Log(\rho_j)$ as the $j-$\textbf{Floquet exponent} of $X(\omega)$. 
\end{definition}
\begin{definition}
    The real parts of Floquet exponents are called \textbf{Lyapunov exponents} and they will be designed as $$\dfrac{1}{\omega}Log(|\rho_j|)=\Re(\lambda_j),\quad j=1,2,\ldots, n.$$
\end{definition}
%%%%%%%%%%%%%%%%%%
As a consequence of the last definition, we have the following result:
\begin{corollary}{\textbf{(Asymptotic behavior of the solutions of a $\omega-$periodic linear IDEPCAG by Floquet exponents})}\label{teo_asympt_behavior_floquet_exponents}

The solutions of \eqref{sistema_idepcag_periodico_homogeneo} converges exponentially to zero if $\Re(\lambda_j)<0$, they will be  $\omega-$periodic if $\Re(\lambda_j)=\Im(\lambda_j)=0$ and
$2\omega-$periodic if $\Re(\lambda_j)=0$ and $\Im(\lambda_j)\neq 0$. Finally, they will be unbounded if $\Re(\lambda_j)>0.$ In other words, if the Lyapunov exponents are less or equal to $0$, solutions of \eqref{sistema_idepcag_periodico_homogeneo} will be bounded. Otherwise, they will be unbounded.  
\end{corollary}

As $X(\omega)$ is non-singular, it has a logarithm. The existence of a logarithm of a matrix is a key fact to establish our version of the Floquet theorem:
\begin{theorem}{\textbf{(Existence of the logarithm of a matrix)}}{(Theorem 2.47)\cite{Chicone}}\label{teo_existencia_logaritmo}
\begin{enumerate}
    \item[1.] If $A$ is a complex nonsingular $n\times n$ matrix, then there exists an $n\times n$ matrix $C$, possibly complex, such that $$\exp{(C)=A}\Leftrightarrow C=Log(A).$$
    \item[2.] If $A$ is a real nonsingular $n\times n$ matrix, then there exists a real $n\times n$ matrix $C$ such that $$\exp{(C)}=A^2 \Leftrightarrow C=Log(A^2).$$ In fact, the real eigenvalues of $A$ will originate positive eigenvalues of $A^2$.
\end{enumerate}
\end{theorem}
%%%%%%%%%%PARA LA TESIS!!
\begin{remark}

    Because it is difficult to find in literature, we will show the importance of the condition $(2)$ of the Last Theorem. Consider the homogeneous linear $\omega-$periodic  ordinary system \eqref{periodico_ordinario}. We see that if all the eigenvalues $\rho_j$ of the Monodromy matrix are real, by the classical Floquet Theorem \ref{teo_floquet_clasico}, we can write a complex solution of \eqref{periodico_ordinario} as $$x_j(t)=\exp{(p_jt)}q(t), \qquad q_j(t)=x_j(t)\exp{(-p_jt)},\quad q_j(t+\omega)=q_j(t),$$ where $p_j=\dfrac{1}{\omega}Log(|\lambda_j|)$ is the monodromy operator and $\rho_j=\exp{(\lambda_j)}$ corresponds to the Floquet multiplier (which is an eigenvalue of the monodromy matrix of \eqref{periodico_ordinario}), i.e. $$x_j(t+\omega)=\rho_j x_j(t).$$ Hence, we have
    \begin{eqnarray*}
        q_j(t+\omega)&=&x_j(t+\omega)\exp{(-(1/\omega)Log(|\lambda_j|)(t+\omega))}\\
        &=&x_j(t)\rho_j \exp{(-(1/\omega)Log(|\lambda_j|)(t+\omega))}\\
        &=&x_j(t)\exp{(Log(\lambda_j)-Log(|\lambda_j|)}\exp{(-(t/\omega)Log(|\lambda_j|))}\\
        &=&sign(\lambda_j)x_j(t)\exp{(-p_jt)}\\
        &=&sign(\lambda_j)q_j(t).
    \end{eqnarray*}
    Now, if we want a real periodic solution of \eqref{periodico_ordinario}, we see that if $\lambda_j\in\mathbb{R}$ is a real eigenvalue of $A$, then we can consider $\Lambda_j=\lambda_j^2$ (i.e., an eigenvalue of $A^2$) to have $\Lambda_j>0$. This way, $q_j(t)$ will be a $2\omega-$periodic function. I.e, if we consider $$\Tilde{p}_j=\dfrac{1}{2\omega}Log(\Lambda_j),$$ then
\begin{eqnarray*}
        \Tilde{q}_j(t+2\omega)&=&x_j(t+2\omega)\exp{(-(1/(2\omega))Log(\Lambda_j)(t+2\omega))}\\
        &=&x_j(t)\rho^2_j \exp{(-(1/(2\omega))Log(|\Lambda_j|)(t+2\omega))}\\
        &=&x_j(t)\exp{(Log(\Lambda_j)-Log(\Lambda_j))}\exp{(-(t/(2\omega))Log(\Lambda_j))}\\
        &=&x_j(t)\exp{(-\Tilde{p}_jt)}\\
        &=&\Tilde{q}_j(t).
\end{eqnarray*}
\end{remark}
By \eqref{monodromia_determinante}, we have the following important result:
\begin{corollary}{\textbf{(L)}}
Let $X(\omega)$ as given in \eqref{MATRIZ_monodromia_avance_retardo}. As
$\det(X(\omega))\neq 0,$  $Log(X(\omega))$ exists.
\end{corollary}
%%%%%%%
Also, if all the related matrices commute, we can give an expression for the logarithm of the monodromy matrix:
\begin{corollary}{\textbf{(LC)}}
Assume that $C_r, A(t),B(t)$ commute for $r=1,\ldots,p$; for all $t\in[0,\omega]$ and $\det(X(\omega))\neq 0.$
%\label{condicion_existencia_log_conmutativo}
Then we have
\begin{eqnarray*}
Log\left(X(\omega)\right)&=&Log(\Phi(\omega,0))+Log\left(
\prod_{r=1}^{p}\left(I+C_{r}\right)\dfrac{J(t_{r},\zeta_{r-1})}{J(t_{r-1},\zeta_{r-1})}
\right).
\label{Logaritmo_MATRIZ_monodromia_avance_retardo}
\end{eqnarray*}
Moreover, for the diagonal case, we have
$$Log\left(X(\omega)\right)=\int_0^{\omega}A(t)dt+Log\left(
\prod_{r=1}^{p}\left(I+C_{r}\right)\dfrac{J(t_{r},\zeta_{r-1})}{J(t_{r-1},\zeta_{r-1})}
\right).$$
\end{corollary}
\begin{proof}
First, as $\det(X(\omega))\neq 0$, we have
\begin{eqnarray*}
Log\left(X(\omega)\right)&=&Log\left(
\prod_{r=1}^{p}\left(I+C_{r}\right) \dfrac{E(t_{r},\zeta_{r-1})}{E(t_{r-1},\zeta_{r-1})}\right).
\label{Logaritmo_MATRIZ_monodromia_avance_retardo_1}
\end{eqnarray*}
Then, as $E(t,s)=\Phi(t,s)J(t,s),$ we see that
\begin{eqnarray*}
Log\left(X(\omega)\right)&=&Log
\left(\prod_{r=1}^{p}\left(I+C_{r}\right) \dfrac{\Phi(t_{r},\zeta_{r-1})}{\Phi(t_{r-1},\zeta_{r-1})}\dfrac{J(t_{r},\zeta_{r-1})}{J(t_{r-1},\zeta_{r-1})}
\right).
\label{Logaritmo_MATRIZ_monodromia_avance_retardo_2}
\end{eqnarray*}
Noting that $\Phi(t_{r},\zeta_{r-1})\phi^{-1}(t_{r-1},\zeta_{r-1})=\Phi(t_{r},t_{r-1})$, we have
\begin{eqnarray*}
Log\left(X(\omega)\right)&=&Log\left(\Phi(\omega,0)\right)+Log\left(
\left(\prod_{r=1}^{p}\left(I+C_{r}\right)\dfrac{J(t_{r},\zeta_{r-1})}{J(t_{r-1},\zeta_{r-1})}\right)
\right).
\label{Logaritmo_MATRIZ_monodromia_avance_retardo_3}
\end{eqnarray*}
Finally, considering the diagonal case, we see that $\Phi(t)=\int_0^t A(u)du$. Hence
\begin{eqnarray}
Log\left(X(\omega)\right)&=&\int_0^\omega A(u)du+Log\left(
\left(\prod_{r=1}^{p}\left(I+C_{r}\right)\dfrac{J(t_{r},\zeta_{r-1})}{J(t_{r-1},\zeta_{r-1})}\right)
\right),
\label{Logaritmo_MATRIZ_monodromia_diagonal}
\end{eqnarray}

and the proof is complete.
\end{proof}
We can now define our $P$ operator:
\begin{equation}
P=\dfrac{1}{\omega}Log\left(X(\omega)\right).
\label{Logaritmo_MATRIZ_monodromia}
\end{equation}
Also, when $C_r, A(t),B(t)$ commute for $r=1,\ldots,p$; for all $t\in[0,\omega]$, we see that 
\begin{eqnarray}
P
%&=&\dfrac{1}{\omega}Log\left(X(\omega)\right)\\
&=&\dfrac{1}{\omega}\left(Log(\Phi(\omega,0))+Log
\left(\prod_{r=1}^{p}\left(I+C_{r}\right)\dfrac{J(t_{r},\zeta_{r-1})}{J(t_{r-1},\zeta_{r-1})}\right)
\right),
\label{Logaritmo_MATRIZ_monodromia_conmutativo}
\end{eqnarray}
and for the diagonal case
\begin{eqnarray}
P=\dfrac{1}{\omega}\left(\int_0^{\omega}A(t)dt+Log
\left(\prod_{r=1}^{p}\left(I+C_{r}\right)\dfrac{J(t_{r},\zeta_{r-1})}{J(t_{r-1},\zeta_{r-1})}\right)
\right),\nonumber
\label{Logaritmo_MATRIZ_monodromia_conmutativo_real}
\end{eqnarray}
where
$\displaystyle{J(t,\tau)=I+\int_{\tau}^{t}\Phi(\tau,s)B(s)ds.}$ 

\begin{remark}
    If $B(t)\coloneqq 0$ and $C_j\coloneqq 0,$ we recover the classical definition of $P$ given in Theorem \ref{teo_floquet_clasico}.
\end{remark}

\section{Main result}
We will state and prove the IDEPCAG version of the Floquet theorem:
\begin{theorem}{\textbf{(Floquet Theorem for IDEPCAG)}}\label{TEO_LYAPUNOV_FLOQUET}

Let the $\omega-$periodic homogeneous linear IDEPCAG \eqref{sistema_idepcag_periodico_homogeneo}: 
\begin{equation*}
\begin{tabular}{ll}
$x^{\prime }(t)=A(t)x(t)+B(t)x(\gamma (t)),$ & $t\neq t_{k},$ \\ 
$\Delta x|_{t=t_{k}}=C_{k}x(t_{k}^{-}),$ & $t=t_{k},$
\end{tabular}
\end{equation*}
and let the conditions \eqref{Periodicidad_coeficientes},\eqref{w-p_propiedad}, Theorem \ref{Teo_forma_normal} and \textbf{(L)} hold. Then, 
\begin{itemize}
\item[(i)]
The solution $X(t)$ of \eqref{sistema_idepcag_periodico_homogeneo} can be represented in the \textbf{Floquet normal form} as
\begin{equation}
X(t)=Q(t)\exp{(Pt)},\quad P=\dfrac{1}{\omega}Log\left(X(\omega)\right),\quad  t\in\mathbb{R},\label{representacion_floquet}
\end{equation}
where $P\in\mathbb{C}^{n\times n}$ is constant and the matrix function  $Q(t)\in \mathcal{PC}^{1}(\mathbb{R},\mathbb{C}^{n\times n})$ is  non-singular, $\omega-$periodic and satisfies the IDEPCAG
\begin{eqnarray}
Q'(t)&=&A(t)Q(t)-Q(t)P+B(t)Q(\gamma(t))e^{P(\gamma(t)-t)}, \qquad \,\,t\neq t_k,\label{q_general}\\
Q(t_k)&=&(I+C_k)Q(t_k^-),\qquad \qquad \qquad \qquad \qquad \qquad \qquad t=t_k.\nonumber
\end{eqnarray}
Also, if $A(t),B(t)$ and $C_k$ are real matrices, each fundamental solution $X(t)$ of \eqref{sistema_idepcag_periodico_homogeneo} can be represented in the Floquet normal form as
\begin{equation}
X(t)=\Tilde{Q}(t)\exp{(\Tilde{P}t)},\quad \Tilde{P}=\dfrac{1}{2\omega}Log(X^2(\omega)),\quad t\in\mathbb{R}, \label{P_y_Q_reales}
\end{equation}
where $\Tilde{P}\in\mathbb{R}^{n\times n}$ is constant and $\Tilde{Q}(t)\in \mathcal{PC}^{1}(\mathbb{R},\mathbb{R}^{n\times n})$ is a non-sigular $2\omega-$periodic matrix function.
\item[(ii)] 
The equation \eqref{sistema_idepcag_periodico_homogeneo} is reducible to the ordinary differential equation:
\begin{equation}
Y'(t)=PY(t),\label{DEPCAG_reducida_floquet}
\end{equation}
by a $\omega-$periodic Floquet-Lyapunov transformation $X(t)=Q(t)Y(t).$ I.e., the IDEPCAG \eqref{sistema_idepcag_periodico_homogeneo} and  \eqref{DEPCAG_reducida_floquet} are \textit{IDEPCAG-Kinematically similar} by the use of the Lyapunov function $Q(t)$, verifying the DEPCAG
$$Q'(t)=A(t)Q(t)-Q(t)P+B(t)Q(\gamma(t))e^{P(\gamma(t)-t)}.$$
\end{itemize}
\end{theorem}
\newpage
\begin{proof}
\begin{itemize}
\item[]
\item[(i)] Since $\det(X(\omega))\neq 0$, by Theorem \ref{teo_existencia_logaritmo} $X(\omega)$ has a logarithm. So, we can rewrite $X(t+\omega)=X(t)X(\omega)$ as $X(t+\omega)=X(t)\exp{(P\omega)},$   with $$P=\dfrac{1}{\omega}Log\left(X(\omega)\right).$$
Now, define 
\begin{equation}
 Q(t)=X(t)\exp{(-Pt)}.\label{forma_definicion_Q}   
\end{equation}
We will prove that the solution of \eqref{sistema_idepcag_periodico_homogeneo} can be written as \eqref{forma_definicion_Q}.\\

First, assuming \eqref{forma_definicion_Q}, we will prove that $Q(t+\omega)=Q(t),\,\forall t\in\mathbb{R}.$\\
Let $X(\omega)\in\mathbb{C}^{n\times n}$ matrix, by Theorem \ref{teo_existencia_logaritmo}, we have
\begin{eqnarray*}
Q(t+\omega)&=&X(t+\omega)\exp{(-P(t+\omega)}\\
&=&X(t)X(\omega)X^{-1}(\omega)\exp{(-Pt)}\\
&=&X(t)\exp{(-Pt)}\\
&=&Q(t).
\end{eqnarray*}
Next, if $X(\omega)\in\mathbb{R}^{n\times n}$, by Theorem \ref{teo_existencia_logaritmo} we define $$\Tilde{P}=\dfrac{1}{2\omega}Log(X(\omega)).$$ Also, we see that $X(t+2\omega)=X(t)X^2(\omega).$ Then,
\begin{eqnarray*}
\Tilde{Q}(t+2\omega)&=&X(t+2\omega)\exp{(-\Tilde{P}(t+2\omega))}\\
&=&X(t)X^2(\omega)(X^{-1}(\omega))^2 \exp{(-\Tilde{P}t)}\\
&=&X(t)\exp{(-\Tilde{P}t)}\\
&=&\Tilde{Q}(t).
\end{eqnarray*}
As $X(t),\exp{(-Pt)}$ and $\exp{(-\Tilde{P}t)}$ are non-singular and differentiable for all $t\in\mathbb{R},$ (possibly with the exceptions at $t=t_k$, when the left-side derivative exists) we have that $Q(t)$ and $\Tilde{Q}(t)$ are non-singular and differentiable too.\\
    
\noindent Now, if we are looking for a solution of the type $X(t)=Q(t)\exp{(Pt)}$ with $Q(t+\omega)=Q(t)$  and $Q(0)=I$, as we will see it has to satisfy \eqref{q_general}.\\
In fact, as $X(t)$ is the solution of \eqref{sistema_idepcag_periodico_homogeneo}, by differentiating the last expression is easy to see that
\begin{eqnarray*}
Q'(t)e^{Pt}+Q(t)Pe^{Pt}&=&A(t)Q(t)e^{Pt}+B(t)Q(\gamma(t))e^{P\gamma(t)},\quad t\neq t_k,\\
\Delta Q(t_k)e^{P(t_k)}&=&C_kQ(t_k^-)e^{(Pt_k)}, \qquad \qquad \qquad \qquad \quad  t=t_k.
\end{eqnarray*}
Multiplying by the right for $\exp{(Pt)}$, we get \eqref{q_general}.\\

Next, following the ideas of \cite{daleckii_krein} (Ch.3), we note that the Cauchy matrix of the solution of the ordinary differential equation $R'(t)=A(t)R(t)-R(t)P$ is $R(t,\tau)=\Phi(t,\tau)R(\tau)\exp{(-P(t-\tau))},$
where $\Phi(t)$ and $\exp{(Pt)}$ are the fundamental matrices of $Z'(t)=A(t)Z(t)$ and $Y(t)=PY(t)$, respectively.\\

For \eqref{q_general}, we have
\begin{equation*}
    Q'(t)-A(t)Q(t)+Q(t)P=B(t)Q(\gamma(t))e^{P(\gamma(t)-t)}.
    %\label{construccionq_1}.
\end{equation*}
Multiplying the last equation for the left by $\Phi(s,t)$, we get
\begin{equation*}
\Phi(s,t)Q'(t)-\Phi(s,t)A(t)Q(t)+\Phi(s,t)Q(t)P=\Phi(s,t)B(t)Q(\gamma(t))e^{P(\gamma(t)-t)}.
\end{equation*}
It is not difficult to see that $\frac{\partial}{\partial t}\Phi(s,t)=-\Phi(s,t)A(t)$. Then, the last equation can be rewritten as
\begin{equation}
\frac{\partial}{\partial t}(\Phi(s,t)Q(t))+\Phi(s,t)Q(t)P=\Phi(s,t)B(t)Q(\gamma(t))e^{P(\gamma(t)-t)}.\label{construccionq_2}
\end{equation}
Next, noting that $\dfrac{d}{dt}(e^{Pt})=Pe^Pt$ and multiplying \eqref{construccionq_2} for the right by $e^{P(t-s)}$, we get
\begin{equation*}
\frac{\partial}{\partial t}(\Phi(s,t)Q(t))e^{P(t-s)}+\Phi(s,t)Q(t)\dfrac{\partial}{\partial t}(e^{P(t-s)})=\Phi(s,t)B(t)Q(\gamma(t))e^{-P(s-\gamma(t))}.
%\label{construccionq_3}
\end{equation*}
%Hence, we can rewrite \eqref{construccionq_3} as
I.e.,
\begin{equation}
    \frac{\partial}{\partial t}(\Phi(s,t)Q(t)e^{P(t-s)})=\Phi(s,t)B(t)Q(\gamma(t))e^{-P(s-\gamma(t))}\notag.
\end{equation}
Now, integrating the last expression from $s$ to $t$, we obtain
\begin{equation*}
    \Phi(s,t)Q(t)e^{P(t-s)}=Q(s)+\displaystyle{\int_s^t}\Phi(s,u)B(u)Q(\gamma(u))e^{-P(s-\gamma(u))}du.
\end{equation*}
Finally, multiplying for the left by $\Phi(t,s)$ and for the right by $e^{-P(t-s)}$ the last equation, we get
\begin{equation}
    Q(t)=\Phi(t,s)Q(s)e^{P(t-s)}+\displaystyle{\int_s^t}\Phi(t,u)B(u)Q(\gamma(u))e^{-P(t-\gamma(u))}du.\label{construccionq_final}
\end{equation}
In the following, we will use \eqref{construccionq_final} rewritten as
\begin{equation}
    Q(t)e^{Pt}=\Phi(t,\tau)Q(\tau)e^{P\tau}+\displaystyle{\int_\tau^t}\Phi(t,u)B(u)Q(\gamma(u))e^{P\gamma(u)}du.\label{q_matriz_fund}
\end{equation}
Using Theorem \ref{TEO_FORMULA_var_PAram} and \eqref{q_matriz_fund}, we will solve \eqref{q_general}.\\
First, let's suppose that $t,\tau\in I_n=[t_n,t_{n+1}),$ for some $n\in\mathbb{Z}.$ In this interval, integrating \eqref{q_general} we get
\begin{equation}
    Q(t)e^{Pt}=\Phi(t,\tau)Q(\tau)e^{P \tau}+\displaystyle{\int_{\tau}^t}\Phi(t,u)B(u)Q(\zeta_n)e^{P\zeta_n}du.\label{q_1}
\end{equation}
Evaluating the last equation at $t=\zeta_n,$ we have
\begin{equation*}
    Q(\zeta_n)e^{P\zeta_n}=\Phi(\zeta_n,\tau)Q(\tau)e^{P \tau}+\displaystyle{\int_{\tau}^{\zeta_n}}\Phi(\zeta_n,u)B(u)Q(\zeta_n)e^{P\zeta_n}du.
    %\label{q_2}
\end{equation*}
We see that
\begin{equation}
    \left(I+\displaystyle{\int_{\zeta_n}^{\tau}}\Phi(\zeta_n,u)B(u)du\right)Q(\zeta_n)e^{P\zeta_n}=\Phi(\zeta_n,\tau)Q(\tau)e^{P \tau}\notag.
\end{equation}
I.e.,
\begin{equation}
Q(\zeta_n)e^{P\zeta_n}=E^{-1}(\tau,\zeta_n)Q(\tau)e^{P\tau}, \label{q_zeta_n}   
\end{equation}
where $E^{-1}(\tau,\zeta_n)=J^{-1}(\tau,\zeta_n)\Phi^{-1}(\tau,\zeta_n).$\\
Now, considering $\tau=\zeta_n$  in \eqref{q_1}, we get
\begin{eqnarray}
    Q(t)e^{Pt}&=&\Phi(t,\zeta_n)Q(\zeta_n)e^{P\zeta_n}+\displaystyle{\int_{\zeta_n}^t}\Phi(t,u)B(u)Q(\zeta_n)e^{P\zeta_n}du\notag\\
    &=&\Phi(t,\zeta_n)\left(I+\displaystyle{\int_{\zeta_n}^t}\Phi(\zeta_n,u)B(u)du\right)Q(\zeta_n)e^{P\zeta_n}\notag\\
    &=&\Phi(t,\zeta_n)J(t,\zeta_n)Q(\zeta_n)e^{P\zeta_n}\notag\\
    &=&E(t,\zeta_n)Q(\zeta_n)e^{P\zeta_n}.\label{q_3}
\end{eqnarray}
Therefore, applying \eqref{q_zeta_n} in \eqref{q_3} we obtain
\begin{equation}
    Q(t)e^{Pt}=E(t,\zeta_n)E^{-1}(\tau,\zeta_n)Q(\tau)e^{P\tau}.\label{q_4}
\end{equation}
Now, evaluating the last equation at $\tau=t_n$, we have
\begin{equation}
     Q(t)e^{Pt}=E(t,\zeta_n)E^{-1}(t_n,\zeta_n)Q(t_n)e^{Pt_n}.\label{q_5}
\end{equation}
Assuming the left-side continuity of the solution, we consider $t\to t_{n+1}^{-}$, getting
\begin{equation*}
     Q(t_{n+1}^{-})e^{P t_{n+1}}=E(t_{n+1}^{-},\zeta_n)E^{-1}(t_n,\zeta_n)Q(t_n)e^{Pt_n}.
     %\label{q_6}
\end{equation*}
Therefore, applying the impulsive condition given by \eqref{q_general}, we get the following difference equation
\begin{equation*}
     Q(t_{n+1})e^{P t_{n+1}}=(I+C_{n+1})E(t_{n+1}^{-},\zeta_n)E^{-1}(t_n,\zeta_n)Q(t_n)e^{Pt_n},
     %\label{q_7}
\end{equation*}
whose solution is
\begin{equation}
Q(t_n)e^{Pt_n}=\left(\displaystyle{\prod_{r=n_0+2}^{n}}(I+C_{r})\dfrac{E(t_{r},\zeta_{r-1})}{E(t_{r-1},\zeta_{r-1})}\right)Q(t_{n_0+1})e^{Pt_{n_0+1}}.\label{q_discreto_casi}    
\end{equation}
By \eqref{q_4}, we see that 
$$Q(t_{n_0+1})e^{Pt_{n_0+1}}=(I+C_{n_0+1})E(t_{n_0+1},\gamma(\tau))E^{-1}(\tau,\gamma(\tau)).$$
So, \eqref{q_discreto_casi} can be rewritten as
\begin{equation}
Q(t_n)e^{Pt_n}=\left(\displaystyle{\prod_{r=n_0+2}^{n}}(I+C_{r})\dfrac{E(t_{r},\zeta_{r-1})}{E(t_{r-1},\zeta_{r-1})}\right)(I+C_{n_0+1})\left(\dfrac{E(t_{n_0+1},\gamma(\tau))}{E(\tau,\gamma(\tau))}\right).\label{q_discreto}    
\end{equation}
Finally, applying \eqref{q_discreto} in \eqref{q_5} we get the solution of \eqref{q_general}:
\begin{eqnarray*}
    Q(t)e^{Pt}&=&\dfrac{E(t,\zeta_{k(t)})}{E(t_{k(t)},\zeta_{k(t)})}\left(\displaystyle{\prod_{r=k(\tau)+2}^{k(t)}}(I+C_{r})\dfrac{E(t_{r},\zeta_{r-1})}{E(t_{r-1},\zeta_{r-1})}\right)\notag\\
    &&\cdot (I+C_{k(\tau)+1})\left(\dfrac{E(t_{k(\tau)+1},\gamma(\tau))}{E(\tau,\gamma(\tau))}\right),
    %\label{q_final}
\end{eqnarray*}
where $k(t)$ is the unique $k\in\mathbb{Z}$ such that $t\in I_{k(t)}=[t_{k(t)},t_{k(t)+1}).$\\

Consequently, as
\begin{eqnarray*}
    X(t)&=&\dfrac{E(t,\zeta_{k(t)})}{E(t_{k(t)},\zeta_{k(t)})}\left(\displaystyle{\prod_{r=k(\tau)+2}^{k(t)}}(I+C_{r})\dfrac{E(t_{r},\zeta_{r-1})}{E(t_{r-1},\zeta_{r-1})}\right)\notag\\
    &&\cdot(I+C_{k(\tau)+1})\left(\dfrac{E(t_{k(\tau)+1},\gamma(\tau))}{E(\tau,\gamma(\tau))}\right),
\end{eqnarray*}
it is straightforward that $$X(t)=Q(t)e^{Pt}.$$
\item[(ii)] 
Finally, by the Floquet-Lyapunov change of variables $X(t)=Q(t)Y(t),$ differentiating at $t\neq t_k$ we have
\begin{eqnarray*}
Q'(t)Y(t)+Q(t)Y'(t)&=&A(t)Q(t)Y(t)+B(t)Q(\gamma(t))Y(\gamma(t))\\
&=&\underbrace{\bigg(Q'(t)+Q(t)P-B(t)Q(\gamma(t))e^{P(\gamma(t)-t)}\bigg)}_{A(t)Q(t) \text{ by } \eqref{q_general}}Y(t)+B(t)Q(\gamma(t))Y(\gamma(t))\\
&=&Q'(t)Y(t)+Q(t)PY(t)-B(t)Q(\gamma(t))e^{P(\gamma(t)-t)}Y(t)\\
&&+B(t)Q(\gamma(t))Y(\gamma(t)).
\end{eqnarray*}
Hence
\begin{eqnarray*}
Q(t)Y'(t)=Q(t)PY(t)-B(t)Q(\gamma(t))e^{P(\gamma(t)-t)}Y(t)+B(t)Q(\gamma(t))Y(\gamma(t)).
\end{eqnarray*}
Since $Q(t)$ is invertible, we have
\begin{equation*}
Y'(t)=PY(t)-Q^{-1}(t)B(t)Q(\gamma(t))e^{P(\gamma(t)-t)}Y(t)+Q^{-1}(t)B(t)Q(\gamma(t))Y(\gamma(t)).
\end{equation*}
\bigskip
Now, for $t=t_k$, by the Floquet normal form, we have
\begin{equation*}
    Q(t_k)e^{(Pt_k)}=(I+C_k)Q(t_k^-)e^{(Pt_k)}.
\end{equation*}
I.e
\begin{equation}
    Q(t_k)=(I+C_k)Q(t_k^-) \label{impulso_Q_exp}
\end{equation}
Also, by the Lyapunov-Floquet change of variables, we have 
\begin{equation*}
    \Delta Q(t_k)Y(t_k)=C_k Q(t_k)Y(t_k^-),
\end{equation*}
i.e.,
\begin{equation*}
    Q(t_k)Y(t_k)=\underbrace{(I+C_k)Q(t_k^-)}_{Q(t_k)}Y(t_k^-).
\end{equation*}
Applying \eqref{impulso_Q_exp} to the last expression and using that $Q(t_k)$ is invertible, we get
$$Y(t_k)=Y(t_k^-).$$
Hence, the impulse effect is not present. So, we reduce the problem to the DEPCAG 
\begin{equation}
Y'(t)=PY(t)-Q^{-1}(t)B(t)Q(\gamma(t))\bigg(e^{P\gamma(t)}e^{-Pt}Y(t)-Y(\gamma(t))\bigg).
\label{DEPCAG_resultante_floquet}
\end{equation}
Now, as $X(t)=Q(t)e^{Pt}$ and $X(t)=Q(t)Y(t),$ then $Y(\gamma(t))=e^{P\gamma(t)}$ and $e^{-Pt}Y(t)=I$. Therefore, rewriting the last equation, 
%\eqref{DEPCAG_resultante_floquet} 
we have
\begin{equation*}
Y'(t)=PY(t).\label{DEPCAG_resultante_floquet_reducida}
\end{equation*}
\end{itemize}
\end{proof}
%%%%%%%%%%%%%%%%%%%%%%%%%%%%%%%%%%%%%%%%%%
\begin{remark}
    It is important to remark that if in \eqref{DEPCAG_resultante_floquet} we consider $\gamma(t)=t,$ then we recover the classical Lyapunov-Floquet equation $$Y'(t)=PY(t).$$
\end{remark}
\begin{corollary}\label{formas_solucion_commutativa_diagonal}
Let Theorem \ref{TEO_LYAPUNOV_FLOQUET} holds. 
\begin{itemize}
%\item[]
\item[(i)] If $A(t),B(t),C_j$ commute $\forall t\in[0,\omega]$ and $j=1,\ldots, p$,  
then $P$ is given by \eqref{Logaritmo_MATRIZ_monodromia_conmutativo_real} and 
\begin{eqnarray*}
    &&P=\dfrac{1}{\omega}\left(Log(\Phi(\omega,0))+Log\left(\displaystyle{\prod_{r=1}^{p}}(I+C_{r})\dfrac{J(t_{r},\zeta_{r-1})}{J(t_{r-1},\zeta_{r-1})}\right)\right),\\
    &&Q(t)=J(t,t_{k(t)})J^{-1}(t_{k(t)},\zeta_{k(t)}),
\end{eqnarray*}
where
$\displaystyle{J(t,\tau)=I+\int_{\tau}^{t}\Phi(\tau,s)B(s)ds.}$ 
\item[(ii)]
If $A(t),B(t),C_j$ are diagonal matrices, then
\begin{eqnarray*}
P=\dfrac{1}{\omega}\left(\displaystyle{\int_0^\omega}A(u)du+\displaystyle{\sum_{r=1}^{p}}Log(\eta_r)\right),\quad Q(t)=\left(\dfrac{I+\displaystyle{\int_{\zeta_{k(t)}}^t} \exp{\left(\displaystyle{\int_s^{\zeta_{k(t)}}A(u)du}\right)}B(s)ds}{I+\displaystyle{\int_{\zeta_{k(t)}}^{t_{k(t)}}}\exp{\left(\displaystyle{\int_s^{\zeta_{k(t)}}A(u)du)}\right)}B(s)ds}\right),
\end{eqnarray*}
where
\begin{equation*}
    \eta_r=(I+C_{r})
\left(\dfrac{I+\displaystyle{\int_{\zeta_{r-1}}^{t_r}} \exp{\left(\displaystyle{\int_s^{\zeta_{r-1}}A(u)du}\right)}B(s)ds}{I+\displaystyle{\int_{\zeta_{r-1}}^{t_{r-1}}}\exp{\left(\displaystyle{\int_s^{\zeta_{r-1}}A(u)du)}\right)}B(s)ds}\right).
\end{equation*}
\item[(iii)]
The Floquet normal form $X(t)=Q(t)\exp{(Pt)}$ of the solution of \eqref{sistema_idepcag_periodico_homogeneo} for the diagonal case is
\begin{eqnarray}
X(t)&=&\left(\dfrac{I+\displaystyle{\int_{\zeta_{k(t)}}^t} \exp{\left(\displaystyle{\int_s^{\zeta_{k(t)}}A(u)du}\right)}B(s)ds}{I+\displaystyle{\int_{\zeta_{k(t)}}^{t_{k(t)}}}\exp{\left(\displaystyle{\int_s^{\zeta_{k(t)}}A(u)du)}\right)}B(s)ds}\right)\exp{\left(\int_0^t A(u)du+ \sum_{r=1}^{k(t)}Log\left(\eta_r\right)\right)}.\notag
\end{eqnarray}   
\end{itemize}
\end{corollary}
\begin{remark}
\begin{itemize}
    \item[] 
    \item If we consider $C_r\coloneqq 0,$ we have the \textbf{DEPCAG} version of Floquet Theory. 
    \item If $C_r\coloneqq B(t)\coloneqq 0,$ we recover the classical version of the Floquet Theorem.
\end{itemize}
\end{remark}
%\begin{equation}
%\exp{(Pt)}=\exp{\left(\int_0^t A(u)du+ \sum_{r=1}^{k(t)}Log\left(\eta_r\right)\right)}.\label{exp_conmutativo}
%\end{equation}
%%%%%%%%%%%%%%%%%%%%%%%%%%%%%%%%%%%%%%%%%%%%%%%%%%%%%%%%%
\begin{corollary}{(Bounded Solution of \eqref{sistema_idepcag_periodico_homogeneo} over $\mathbb{R}$)}

    The only bounded solution of \eqref{sistema_idepcag_periodico_homogeneo} over all $\mathbb{R}$ is the $\omega-$periodic or the $2\omega-$periodic solution. I.e., when the Lyapunov exponent is $0,$ but the Floquet exponent is purely imaginary or when the Floquet exponent is identically $0$.
\end{corollary}

 \begin{remark}
The problem of finding a normal form of the Floquet solution of \eqref{sistema_idepcag_periodico_homogeneo} is equivalent to finding $P$ and $Q(t)$ satisfying \eqref{q_general}. In general, this problem seems to be very difficult. (See \cite{Castelli_lessard}).\\
 \end{remark}
%%%%%%%%%%%%

\section{Some examples}
Let the following $1-$periodic IDEPCA
\begin{equation}
\begin{tabular}{ll}
$z'(t)=\sin(2\pi t)z\left(\left[t\right]\right),$ & $t\neq k,\quad k\in\mathbb{N}$, \\ 
$z(k)=c z(k^{-}),$ & $t=k,$ \\ 
$z(0)=1$. & 
\end{tabular}
\label{IDEPCA_EJEMPLO}
\end{equation}
We see that $\sigma_k^{-}(0)=\nu_k^{-}(\sin(2\pi t))=0<1,$ and $
J(t,\tau)=E(t,\tau)=1+\int_{\tau}^{t}\sin(2\pi s)ds.$\\
As $\displaystyle{\int_{j}^{j+1}\sin\left(2\pi s\right)ds=0,\,\, \forall j\in\mathbb{Z},}$
by Corollary \ref{formas_solucion_commutativa_diagonal}, the solution of \eqref{IDEPCA_EJEMPLO} is
%\begin{equation*}
%z(t)=
%c^{[t]}\left(1+\int_{[t]}^{t}\sin\left(2\pi s\right)ds\right).
%\label{MATRIZ_FUNDAMENTAL_IDEPCA_ejemplo}
%\end{equation*}
%Hence, the solution of \eqref{IDEPCA_EJEMPLO} with $z(0)=1$ is
$z(t)=c^{[t]}\left(1+\int_{[t]}^{t}\sin(2\pi s)ds\right),$
or
$$z(t)=\exp{\left(Log\left(c\right)[t]\right )}\left(1+ \dfrac{\cos(2 \pi [t]) - \cos(2 \pi t)}{2 \pi}\right).$$
Moreover, by Corollaries \ref{teo_asympt_behavior_floquet_multipliers} and \ref{teo_asympt_behavior_floquet_exponents}, we have the following description of the asymptotic behavior of the solutions:
%\newpage
\begin{enumerate}
\item[(i)] if $c=-4/5$, the Lyapunov exponent of the system is $ln(4/5)<0.$ So, the zero solution is exponentially asymptotically stable.
\item[(ii)] if $c=1.1$, the Lyapunov exponent of the system is $ln(1.1)>0$. Consequently, the solution is unbounded.
\item[(iii)] if $c=-1$, the Floquet multiplier satisfies $|\rho|=1$, and the Lyapunov exponent is $0$, but the imaginary part of the Floquet exponent is non-zero. Therefore, the solution is $2-$periodic and oscillatory. We remark that if $\Im(\lambda)\neq 0$, then there is an oscillatory solution.
\item[(iv)] if $c=1$ (non-impulsive case), the Floquet multiplier satisfies $|\rho|=1$, and the Lyapunov exponent is $0$. In this case, the Floquet exponent is equal to $0$. Hence, the solution is $1-$periodic.
\end{enumerate}
%%%%%FIGURAS
\begin{figure}[h!]
\centering
\includegraphics[scale=0.23]{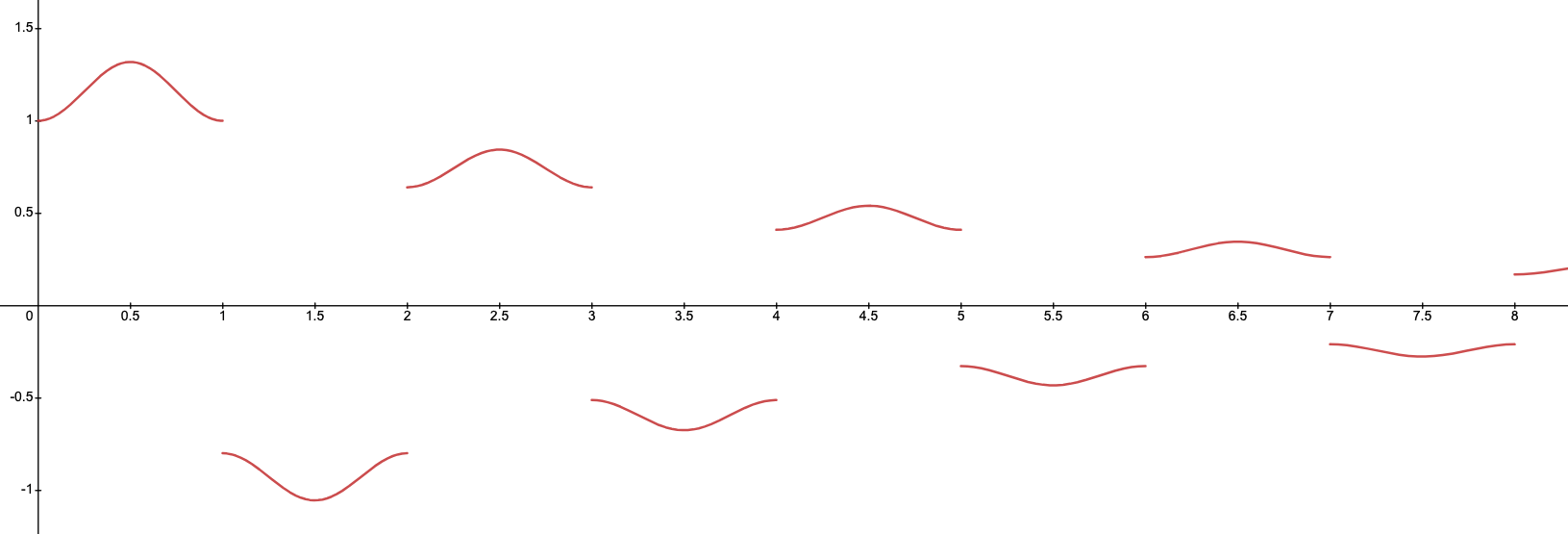}
\caption{Solution of  \eqref{IDEPCA_EJEMPLO} with $c=-4/5$ and $z_0=1$.}
\end{figure}
%%%%%
\begin{figure}[h!]
\centering
\includegraphics[scale=0.39]{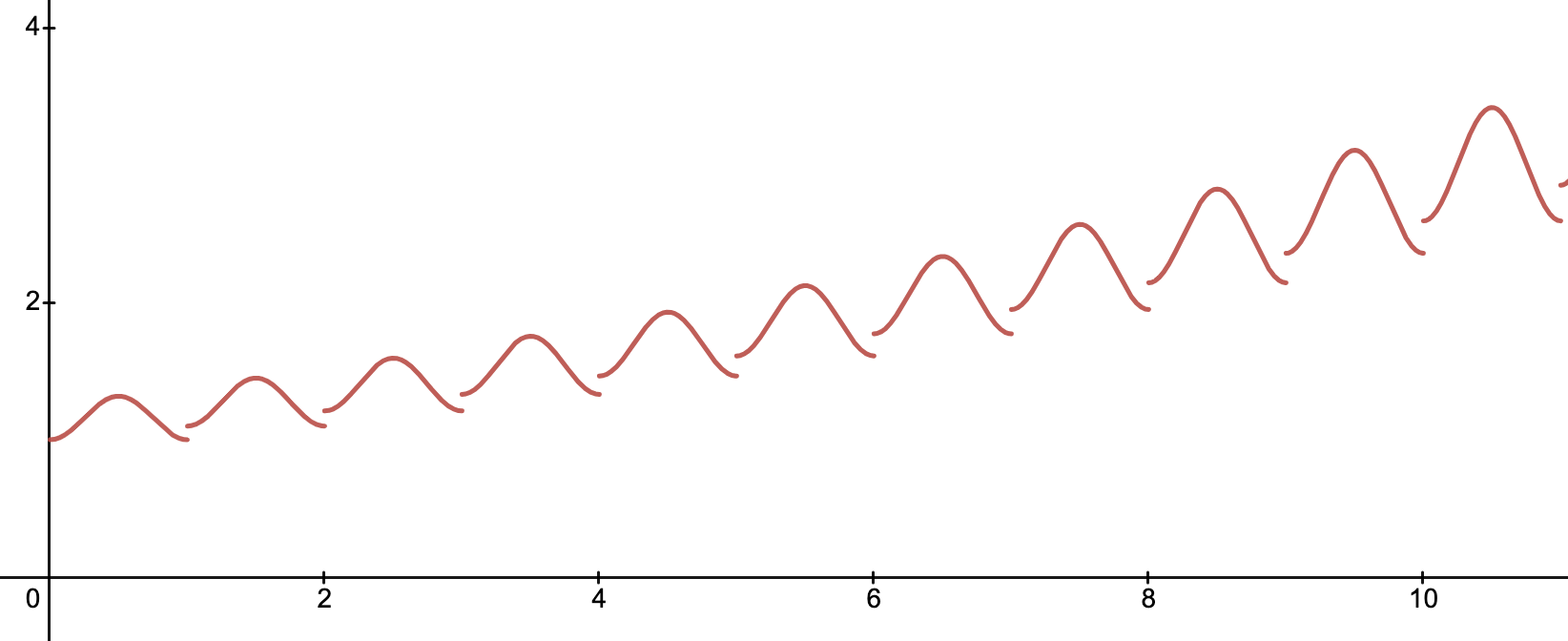}
\caption{Solution of  \eqref{IDEPCA_EJEMPLO} with $c=1.1$ and $z_0=1$.}
\end{figure}
%%%%%%%%
\begin{figure}[h!]
\centering
\includegraphics[scale=0.3]{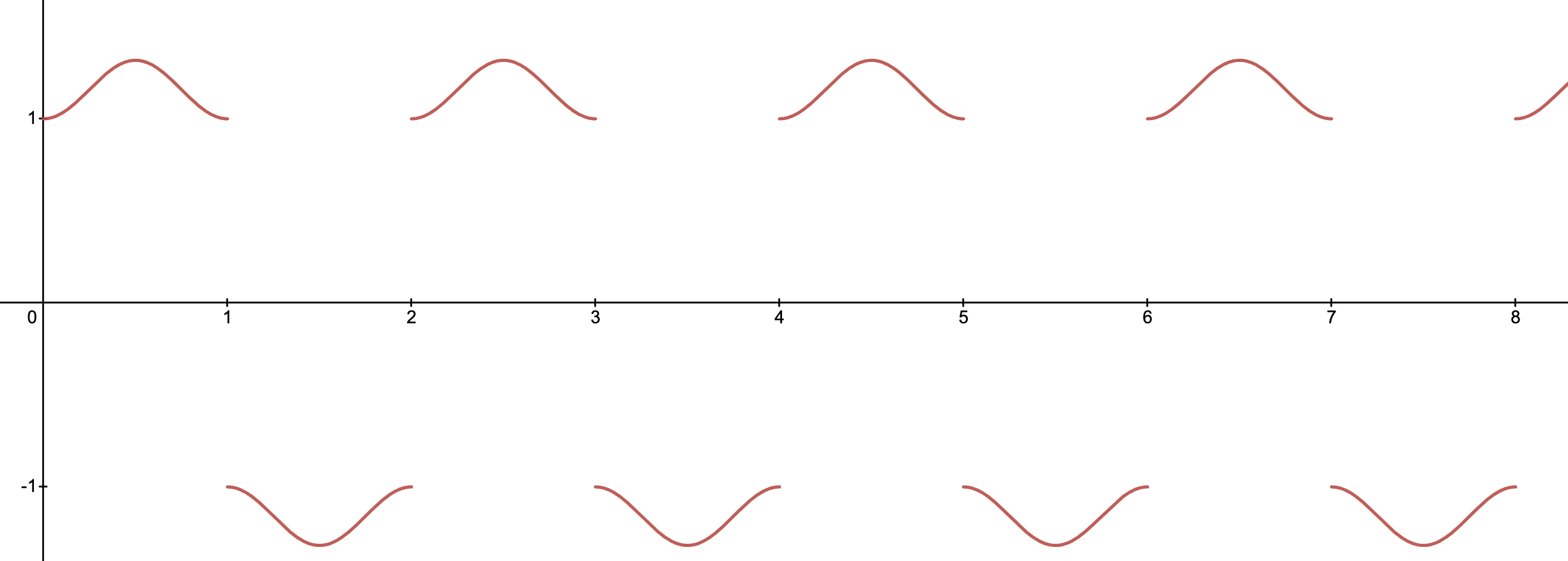}
\caption{Solution of  \eqref{IDEPCA_EJEMPLO} with $c=-1$ and $z_0=1$.}
\end{figure}
%%%%%%%%%%%
\begin{figure}[h!]
\centering
\includegraphics[scale=0.39]{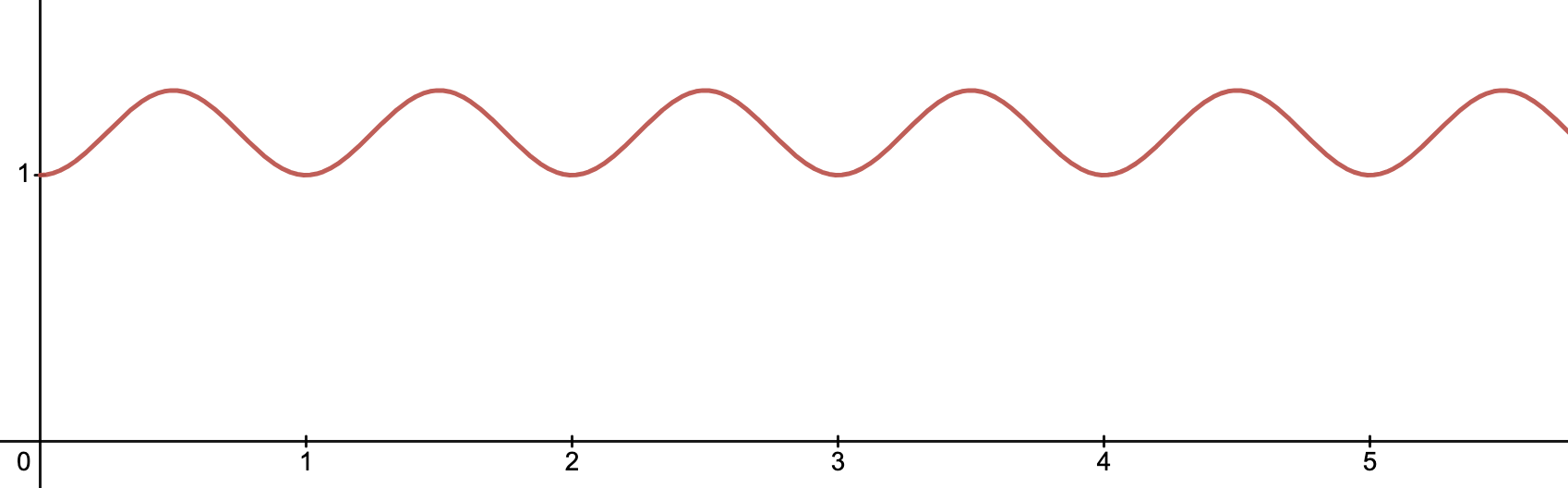}
\caption{Solution of  \eqref{IDEPCA_EJEMPLO} with $c=1$ and $z_0=1$.}
\end{figure}
\newpage

%%%%%%%%%%%%%%%%%%%%%%%%%%%%%%%%%%%%%%%%%
\begin{example}
Inspired in Ex. 3.2 of \cite{Folkers}, 
let the following IDEPCA system
\begin{equation}
\begin{tabular}{ll}
$X'(t)=A(t)X(t)+B(t)X(2\pi\left[\frac{t}{2\pi}\right]),$ & $t\neq 2k\pi,\quad k\in\mathbb{Z}$, \\ 
$X(2k\pi)=C X(2k\pi^{-}),$ & $t=2k\pi,$ \\ 
$X(0)=I$, & 
\end{tabular}
\label{IDEPCA_EJEMPLO_matricial}
\end{equation}
where
$$A(t)=
\begin{pmatrix}
\cos(t) & -\sin(t)\\
&\\
\sin(t) & \cos(t)
\end{pmatrix},
\quad 
B(t)=
\begin{pmatrix}
1 & 0\\
&\\
0 & 1
\end{pmatrix}
\text{   and   } 
C=
\begin{pmatrix}
\frac{1}{5} & 0\\
&\\
0 & \frac{1}{5}
\end{pmatrix}
.$$
The matrix $A(t)$ is $2\pi-$periodic and $\gamma(t)=2\pi[t/2\pi]$ verifies
\begin{eqnarray*}
\gamma(t)=2k\pi, \text{ when } t\in I_k=[2k\pi,2(k+1)\pi),\,\,k\in\mathbb{Z}.
\end{eqnarray*}
Hence, we have $t_k=\zeta_k=2k\pi$ and
\begin{eqnarray*}
    t_{k+1}=t_k+2\pi,\label{2pi_1_propiedad}
    \quad \zeta_{k+1}=\zeta_k+2\pi.
\end{eqnarray*}
The ordinary system $Z'(t)=A(t)Z(t)$ has
\begin{eqnarray*}
\Phi(t)=&&
\begin{pmatrix}
e^{\sin(t)}\cos(1-\cos(t)) & -e^{\sin(t)}\sin(1-\cos(t))\\
&\\
e^{\sin(t)}\sin(1-\cos(t)) & e^{\sin(t)}\cos(1-\cos(t))
\end{pmatrix}
\\
=&&\underbrace{\begin{pmatrix}
e^{\sin(t)} & 0\\
&\\
0 & e^{\sin(t)}
\end{pmatrix}}_{M(t)}
\underbrace{\begin{pmatrix}
\cos(1-\cos(t)) & -\sin(1-\cos(t))\\
&\\
\sin(1-\cos(t)) & \cos(1-\cos(t))
\end{pmatrix}}_{N(t)},
\end{eqnarray*}
as the fundamental matrix satisfying $\Phi(0)=I.$\\
Also, as $M(t)N(t)=N(t)M(t)$, we see that 
\begin{eqnarray*}
Log(\Phi(2\pi))&=&Log(M(2\pi))+Log(N(2\pi))\\
&=&2Log(I)\\
&=&
\begin{pmatrix}
0 & 0\\
&\\
0 & 0
\end{pmatrix}.
\end{eqnarray*}
As
$$\Phi^{-1}(t)=
\begin{pmatrix}
\exp{(-\sin(t))}(\cos(1-\cos(t))) & \exp{(-\sin(t))}(\sin(1-\cos(t)))\\
&\\
-\exp{(-\sin(t))}(\sin(1-\cos(t))) & \exp{(-\sin(t))}(\cos(1-\cos(t)))
\end{pmatrix},
$$
$B(t)=I$ and  $J(0,2\pi)=I+\int_{0}^{2\pi} \Phi^{-1}(s)ds,$ by Corollary \ref{formas_solucion_commutativa_diagonal}, we have
\begin{eqnarray*}
X(2\pi)&&=C\bigg(I+\displaystyle{\int_0^{2\pi}}\Phi^{-1}(s)ds\bigg)\\
&&=
\begin{pmatrix}
i & -i\\
&\\
1 & 1
\end{pmatrix}
\begin{pmatrix}
0.878964-1.05742i & 0\\
&\\
0 & 0.878964+1.05742i
\end{pmatrix}
\begin{pmatrix}
-\frac{i}{2} & \frac{1}{2}\\
&\\
\frac{i}{2} & \frac{1}{2}
\end{pmatrix}.
\end{eqnarray*}
In this way, we have
\begin{eqnarray*}
Log\left(X(2\pi)\right)=
\begin{pmatrix}
i & -i\\
&\\
1 & 1
\end{pmatrix}
\begin{pmatrix}
Log(0.878964-1.05742i) & 0\\
&\\
0 & Log(0.878964+1.05742i)
\end{pmatrix}
\begin{pmatrix}
-\frac{i}{2} & \frac{1}{2}\\
&\\
\frac{i}{2} & \frac{1}{2}
\end{pmatrix}.
\end{eqnarray*}
Hence, we get
$$
P=\dfrac{1}{2\pi}Log\left(X(2\pi)\right)=
\begin{pmatrix}
0.0253436-0.0698132i & 0\\
&\\
0 & 0.0253436+0.0698132i
\end{pmatrix}
$$
Therefore, as $0.0253436>0$, by corollary \ref{teo_asympt_behavior_floquet_exponents} the solutions of system \eqref{IDEPCA_EJEMPLO_matricial} are \textbf{unbounded}.\\
Finally, by Corollary \ref{formas_solucion_commutativa_diagonal}, the Floquet normal form of the solutions of \eqref{IDEPCA_EJEMPLO_matricial} is $X(t)=Q(t)e^{Pt},$
where 
$$Q(t)=
\begin{pmatrix}
1+\displaystyle{\int_{2\pi\left[\frac{t}{2\pi}\right]}^t} \exp{(-\sin(s))}(\cos(1-\cos(s)))ds & \displaystyle{\int_{2\pi\left[\frac{t}{2\pi}\right]}^t} \exp{(-\sin(s))}(\sin(1-\cos(s)))ds\\
&\\
-\displaystyle{\int_{2\pi\left[\frac{t}{2\pi}\right]}^t} \exp{(-\sin(s))}(\sin(1-\cos(s)))ds & 1+\displaystyle{\int_{2\pi\left[\frac{t}{2\pi}\right]}^t} \exp{(-\sin(s))}(\cos(1-\cos(s)))ds
\end{pmatrix},$$
and
$$
e^{Pt}=
\begin{pmatrix}
\exp{((1.02317-0.0715469i)t)} & 0\\
&\\
0 & \exp{((1.02317+0.0715469i)t)}
\end{pmatrix}.
$$
If we consider 
$X(t)=(x_1(t)\,\,x_2(t))^{t}
$
with $X(0)=(0\,\,1)^{t},
$
 the solution of \eqref{IDEPCA_EJEMPLO_matricial} is
\begin{eqnarray*}
\begin{pmatrix}
x_1(t)\\
&\\
x_2(t)
\end{pmatrix}=
\begin{pmatrix}
\exp{((1.02317+0.0715469i)t)}\left(\displaystyle{\int_{2\pi\left[\frac{t}{2\pi}\right]}^t} \exp{(-\sin(s))}(\sin(1-\cos(s)))ds\right)\\
&\\
\exp{((1.02317+0.0715469i)t)}\left(\displaystyle{1+\int_{2\pi\left[\frac{t}{2\pi}\right]}^t} \exp{(-\sin(s))}(\cos(1-\cos(s)))ds \right)
\end{pmatrix}
\end{eqnarray*}
which is clearly unbounded.
\end{example}
%%%%%%%%%%%%%%%% 3d real %%%%%%%%%%
\begin{figure}[h!]
\centering
\includegraphics[scale=0.42]{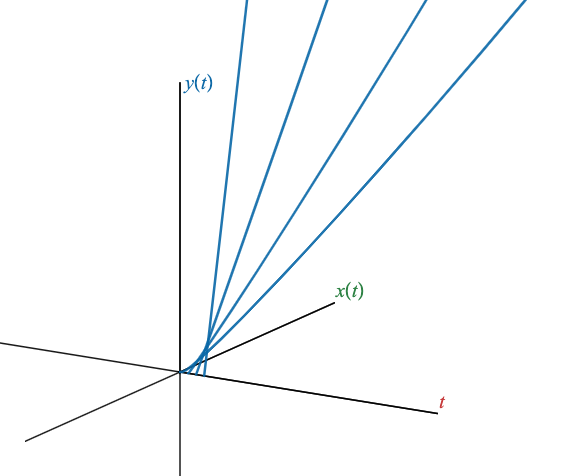}
\caption{Solution of  \eqref{IDEPCA_EJEMPLO_matricial} with  $X_0=(0,1)^{t},\,f(t)=(t,\Re(x_1(t)),\Re(x_2(t)))$.}
\end{figure}
%%%%%%%%%%%%%%% 3d  imaginario
\begin{figure}[h!]
\centering
\includegraphics[scale=0.45]{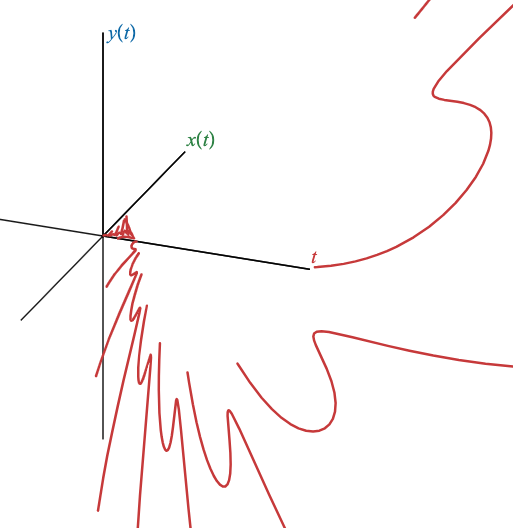}
\caption{Solution of  \eqref{IDEPCA_EJEMPLO_matricial} with  $X_0=(0,1)^{t},\,g(t)=(t,\Im(x_1(t)),\Im(x_2(t)))$.}
\end{figure}
%%%%%%%%%%%%%%%  2d  x1 %%%%%%%%%%%%%%%%%%%%%
%\begin{figure}[h!]
%\centering
%\includegraphics[scale=0.45]{2d_Floquet_real_x1.png}
%\caption{Solution of  \eqref{IDEPCA_EJEMPLO_matricial} with  $X_0=(0,1)^{t},\,f(t)=(t,\Re(x(t)),\Re(y(t)))$.}
%\end{figure}
%%%%%
%\begin{figure}[h!]
%\centering
%\includegraphics[scale=0.4]{2d_Floquet_im_x1.png}
%\caption{Solution of  %\eqref{IDEPCA_EJEMPLO_matricial} with  $X_0=(0,1)^{t},\,g(t)=\Im(x_1(t))$.}
%\end{figure}

%%%%%%%%%%%%%%%  2d  x2 %%%%%%%%%%%%%%%%%%%%%
%\begin{figure}[h!]
%\centering
%\includegraphics[scale=0.45]{2d_Floquet_real_x2.png}
%\caption{Solution of  \eqref{IDEPCA_EJEMPLO_matricial} with  $X_0=(0,1)^{t},\,f(t)=\Re(x_2(t)).$}
%\end{figure}
%%%%%
%\begin{figure}[h!]
%\centering
%\includegraphics[scale=0.4]{2d_Floquet_im_x2.png}
%\caption{Solution of  \eqref{IDEPCA_EJEMPLO_matricial} with  $X_0=(0,1)^{t},\,g(t)=\Im(x_2(t))$.}
%\end{figure}
\newpage
%%%%%%%%
%%%%%%%%%5
\section{Conclusions}
Our research presented a version of the classical Lyapunov-Floquet Theorem for nonautonomous linear impulsive differential equations with piecewise constant arguments of generalized type. To the best of our knowledge, this marks the first extension of the Floquet Theory to this particular class of differential equations.

\section*{Acknowledgments}
\emph{Ricardo Torres} sincerely thanks Prof. Manuel Pinto for the encouragement to work in this subject and for all his support during my career. Also, the author sincerely thanks DESMOS PBC for granting permission to use the images employed in this work. They were created with the DESMOS graphic calculator \\ \url{https://www.desmos.com/calculator}.
%\section*{Funding}
%This research received no specific grant from funding agencies in the public, commercial, or not-for-profit sectors.

\end{document}